\documentclass[10pt,a4paper]{article}

\usepackage[leqno]{amsmath}
\usepackage{graphicx}
\usepackage{latexsym}
\usepackage{amsmath,amsfonts,amssymb,amsthm,mathrsfs,euscript,makeidx,color}
\usepackage{enumerate}

\usepackage[colorlinks,linkcolor=blue,anchorcolor=green,citecolor=red]{hyperref}

\def\sqr#1#2{{\vcenter{\vbox{\hrule height.#2pt
              \hbox{\vrule width.#2pt height#1pt \kern#1pt \vrule width.#2pt}
              \hrule height.#2pt}}}}
\def\signed #1{{\unskip\nobreak\hfil\penalty50
              \hskip2em\hbox{}\nobreak\hfil#1
              \parfillskip=0pt \finalhyphendemerits=0 \par}}
\def\endpf{\signed {$\sqr69$}}

\def\5n{\negthinspace \negthinspace \negthinspace \negthinspace \negthinspace }
\def\4n{\negthinspace \negthinspace \negthinspace \negthinspace }
\def\3n{\negthinspace \negthinspace \negthinspace }
\def\2n{\negthinspace \negthinspace }
\def\1n{\negthinspace }

\def\dbE{\mathbb{E}}
\def\dbF{\mathbb{F}}

\def\dbH{\mathbb{H}}

\def\dbP{\mathbb{P}}

\def\dbR{\mathbb{R}}
\def\dbS{\mathbb{S}}

\def\dbX{\mathbb{X}}
\def\dbY{\mathbb{Y}}
\def\dbZ{\mathbb{Z}}

\def\sR{\mathscr{R}}

\def\cF{{\cal F}}

\def\cQ{{\cal Q}}

\def\cU{{\cal U}}

\def\BA{{\bf A}}
\def\BB{{\bf B}}
\def\BC{{\bf C}}
\def\BD{{\bf D}}

\def\BG{{\bf G}}

\def\BI{{\bf I}}
\def\BJ{{\bf J}}

\def\BP{{\bf P}}
\def\BQ{{\bf Q}}
\def\BR{{\bf R}}
\def\BS{{\bf S}}

\def\BY{{\bf Y}}
\def\BZ{{\bf Z}}

\def\Bg{{\bf g}}

\def\Bq{{\bf q}}

\def\Brho{{\rho\3n\rho\3n\rho\3n\rho\3n\rho}}

\def\BTh{\,{\Th\3n\2n\Th\3n\2n\Th}}
\def\BBeta{{\eta\3n\eta\3n\eta}}
\def\BBzeta{{\z\3n\z\3n\z\3n\z\3n\z\3n\z\3n\z}}
\def\BPi{{\bf\Pi}}

\def\ae{\hbox{\rm a.e.}}
\def\as{\hbox{\rm a.s.}}

\def\ds{\displaystyle}

\def\ns{\noalign{\ss}}
\def\no{\noindent}

\def\ss{\smallskip}
\def\ms{\medskip}

\def\q{\quad}
\def\qq{\qquad}
\def\hb{\hbox}

\def\sc{\scriptstyle}

\def\({\Big (}
\def\){\Big )}
\def\[{\Big[}
\def\]{\Big]}
\def\lt{\left}
\def\rt{\right}
\def\lan{\langle}
\def\ran{\rangle}

\def\blan{\big\langle}
\def\bran{\big\rangle}

\def\rf{\eqref}

\def\a{\alpha}
\def\b{\beta}

\def\e{\varepsilon}
\def\z{\zeta}

\def\l{\lambda}

\def\n{\nu}
\def\si{\sigma}

\def\f{\varphi}
\def\th{\theta}

\def\i{\infty}
\def\G{\Gamma}

\def\Th{\Theta}

\def\O{\Omega}

\def\Ra{\mathop{\Rightarrow}}

\def\esssup{\mathop{\rm esssup}}

\def\wt{\widetilde}

\def\cd{\cdot}

\def\tr{\hbox{\rm tr$\,$}}

\def\deq{\triangleq}
\def\les{\leqslant}
\def\ges{\geqslant}
\def\bde{\begin{definition}\label}
\def\ede{\end{definition}}
\def\be{\begin{equation}}
\def\bel{\begin{equation}\label}
\def\ee{\end{equation}}
\def\bt{\begin{theorem}\label}
\def\et{\end{theorem}}
\def\bc{\begin{corollary}\label}
\def\ec{\end{corollary}}
\def\bl{\begin{lemma}\label}
\def\el{\end{lemma}}
\def\bp{\begin{proposition}\label}
\def\ep{\end{proposition}}
\def\bas{\begin{assumption}\label}
\def\eas{\end{assumption}}
\def\br{\begin{remark}\label}
\def\er{\end{remark}}
\def\bex{\begin{example}\label}
\def\ex{\end{example}}
\def\ba{\begin{array}}
\def\ea{\end{array}}
\def\ben{\begin{enumerate}}
\def\een{\end{enumerate}}

\def\square#1{\vbox{\hrule\hbox{\vrule height#1%
     \kern#1\vrule}\hrule}}
\def\rectangle#1#2{\vbox{\hrule\hbox{\vrule height#1%
     \kern#2\vrule}\hrule}}

\font\tenbb=msbm10 \font\sevenbb=msbm7 \font\fivebb=msbm5

\newfam\bbfam
\scriptscriptfont\bbfam=\fivebb \textfont\bbfam=\tenbb
\scriptfont\bbfam=\sevenbb

\newtheorem{theorem}{\indent Theorem}[section]
\newtheorem{definition}[theorem]{\indent Definition}
\newtheorem{proposition}[theorem]{\indent Proposition}
\newtheorem{corollary}[theorem]{\indent Corollary}
\newtheorem{lemma}[theorem]{\indent Lemma}
\newtheorem{remark}[theorem]{\indent Remark}
\newtheorem{example}[theorem]{\indent Example}

\newtheorem{assumption}[theorem]{\indent Assumption}

\makeatletter
   
   \@addtoreset{equation}{section}
\makeatother

\begin{document}

\title{\bf Linear Quadratic Stochastic Two-Person\\ Nonzero-Sum Differential Games:\\
Open-Loop and Closed-Loop Nash Equilibria\thanks{This work is supported in part by NSF Grant
DMS-1406776.}}

\author{Jingrui Sun\thanks{Department of Applied Mathematics, The Hong Kong Polytechnic University, Hong Kong, China (sjr@mail.ustc.edu.cn).}\,\, \ \ and \ \
Jiongmin Yong\thanks{Department of Mathematics, University of Central Florida, Orlando, FL 32816, USA (jiongmin.yong@ucf.edu).}}
\maketitle

\no\bf Abstract: \rm
In this paper, we consider a linear quadratic stochastic two-person nonzero-sum differential game.
Open-loop and closed-loop Nash equilibria are introduced. The existence of the former is characterized
by the solvability of a system of forward-backward stochastic differential equations, and that of the
latter is characterized by the solvability of a system of coupled symmetric Riccati differential equations.
Sometimes, open-loop Nash equilibria admit a closed-loop representation, via the solution to a system of non-symmetric Riccati equations, which is different from the outcome of the closed-loop Nash equilibria in general. However, it is found that for the case of zero-sum differential games, the Riccati equation system for the closed-loop representation of open-loop saddle points coincides with that for the closed-loop saddle points, which leads to the conclusion that the closed-loop representation of open-loop saddle points is the outcome of the corresponding closed-loop saddle point as long as both exist. In particular, for linear quadratic optimal control problem, the closed-loop representation of open-loop optimal controls coincides with the outcome of the corresponding closed-loop optimal strategy, provided both exist.

\ms

\no\bf Keywords: \rm stochastic differential equation, linear quadratic differential game, two-person,
nonzero-sum, Nash equilibrium, Riccati differential equation, closed-loop, open-loop.

\ms

\no\bf AMS Mathematics Subject Classification. \rm 93E20, 91A23, 49N70, 49N10.

\section{Introduction}

Let $(\O,\cF,\dbF,\dbP)$ be a complete filtered probability space on which a standard one-dimensional
Brownian motion $\{W(t),t\ges0\}$ is defined such that $\dbF=\{\cF_t\}_{t\ges0}$ is the natural
filtration of $W(\cd)$ augmented by all the $\dbP$-null sets in $\cF$.
Consider the following controlled linear (forward) stochastic differential equation (FSDE, for short)
on $[t,T]$:
\bel{state}\left\{\2n\ba{ll}
\ds dX(s)=\big[A(s)X(s)+B_1(s)u_1(s)+B_2(s)u_2(s)+b(s)\big]ds\\
\ns\ds\qq\q~~~+\big[C(s)X(s)+D_1(s)u_1(s)+D_2(s)u_2(s)+\si(s)\big]dW(s),\qq s\in[t,T],\\
\ns\ds X(t)= x.\ea\right.\ee
In the above, $X(\cd)$ is called the {\it state process} taking values in the $n$-dimensional
Euclidean space $\dbR^n$ with the {\it initial pair} $(t,x)\in[0,T)\times\dbR^n$; for $i=1,2$,
$u_i(\cd)$ is called the {\it control process} of Player $i$ taking values in $\dbR^{m_i}$. We
assume that the {\it coefficients} $A(\cd)$, $B_1(\cd)$, $B_2(\cd)$, $C(\cd)$, $D_1(\cd)$, and
$D_2(\cd)$ are deterministic matrix-valued functions of proper dimensions, and that $b(\cd)$
and $\si(\cd)$ are $\dbF$-progressively measurable processes taking values in $\dbR^n$.
For $i=1,2$ and $t\in[0,T)$, we define
$$\cU_i[t,T]=\Big\{u_i:[t,T]\times\O\to\dbR^{m_i}\bigm|u_i(\cd)
\hb{ is $\dbF$-progressively measurable, }\dbE\int_t^T|u_i(s)|^2ds<\i\Big\}.$$
Any element $u_i(\cd)\in\cU_i[t,T]$ is called an {\it admissible control} of Player $i$ on $[t,T]$. Under some mild conditions on the coefficients, for any initial pair $(t,x)\in[0,T)\times\dbR^n$
and controls $u_i(\cd)\in\cU_i[t,T]$, $i=1,2$, the state equation \rf{state} admits a unique
solution $X(\cd)\equiv X(\cd\,;t,x,u_1(\cd),u_2(\cd))$.
The cost functional for Player $i$ is defined by the following:
\bel{cost}\ba{ll}
\ds J^i(t,x;u_1(\cd),u_2(\cd))\deq\dbE\Big\{\lan G^iX(T),X(T)\ran+2\lan g^i,X(T)\ran\\
\ns\ds\qq\qq\qq\qq\q+\int_t^T\[\lan
{\sc\begin{pmatrix}\sc Q^i(s)   &\1n\sc S^i_1(s)^\top &\1n\sc S^i_2(s)^\top\\
                   \sc S^i_1(s) &\1n\sc R^i_{11}(s)   &\1n\sc R^i_{12}(s)  \\
                   \sc S^i_2(s) &\1n\sc R^i_{21}(s)   &\1n\sc R^i_{22}(s)  \end{pmatrix}}
{\sc\begin{pmatrix}\sc X(s) \\ \sc u_1(s) \\ \sc u_2(s)\end{pmatrix}},
{\sc\begin{pmatrix}\sc X(s) \\ \sc u_1(s) \\ \sc u_2(s)\end{pmatrix}}\ran
+2\lan
{\sc\begin{pmatrix}\sc q^i(s) \\ \sc\rho^i_1(s) \\ \sc\rho^i_2(s)\end{pmatrix},
    \begin{pmatrix}\sc X(s)   \\ \sc u_1(s)     \\ \sc u_2(s)    \end{pmatrix}}\ran\]ds\Big\},\ea\ee
where $Q^i(\cd)$, $S^i_1(\cd)$, $S^i_2(\cd)$, $R^i_{11}(\cd)$, $R^i_{12}(\cd)$, $R^i_{21}(\cd)$, and $R^i_{22}(\cd)$ are deterministic matrix-valued functions of proper dimensions with
$$Q^i(\cd)^\top=Q^i(\cd),\q R^i_{jj}(\cd)^\top=R^i_{jj}(\cd),
\q R^i_{12}(\cd)^\top=R^i_{21}(\cd),  \qq i,j=1,2,$$
where the superscript $^\top$ denotes the transpose of matrices, and $G^i$ is a symmetric matrix;
$q^i(\cd)$, $\rho^i_1(\cd)$, and $\rho^i_2(\cd)$ are allowed
to be vector-valued $\dbF$-progressively measurable processes, and $g^i$ is allowed to be an
$\cF_T$-measurable random vector. Then we can formally pose the following problem.

\ms

\bf Problem (SDG). \rm For any initial pair $(t,x)\in[0,T)\times\dbR^n$ and $i=1,2$,
Player $i$ wants to find a control $u^*_i(\cd)\in\cU_i[t,T]$ such that the cost
functional $J^i(t,x;u_1(\cd),u_2(\cd))$ is minimized.

\ms

The above posed problem is referred to as a linear quadratic (LQ, for short) stochastic
{\it two-person differential game}. In the case
\bel{16Apr2-zero}\ba{ll}
\ds J^1(t,x;u_1(\cd),u_2(\cd))+J^2(t,x;u_1(\cd),u_2(\cd))=0,\\
\ns\ds\qq\qq\qq\qq\qq\qq\forall\,(t,x)\in[0,T]\times\dbR^n,~\forall\,u_i(\cd)\in\cU_i[t,T],\q i=1,2,
\ea\ee
the corresponding Problem (SDG) is called an LQ stochastic two-person {\it zero-sum}
differential game. To guarantee \rf{16Apr2-zero}, one usually assumes that
\bel{16Apr2-1+2=0}\ba{lllll}
\ds G^1+G^2=0, &~ g^1+g^2=0, &~ Q^1(\cd)+Q^2(\cd)=0, &~ q^1(\cd)+q^2(\cd)=0,\\
\ns\ds S_j^1(\cd)+S^2_j(\cd)=0,       &~ R^1_{jk}(\cd)+R^2_{jk}(\cd)=0,
    &~ \rho^1_j(\cd)+\rho^2_j(\cd)=0, &~ j,k=1,2.
\ea\ee
We refer the readers to \cite{Sun-Yong 2014} (and the references cited therein) for the case of
LQ stochastic two-person zero-sum differential games. Recall that in \cite{Sun-Yong 2014},
open-loop and closed-loop saddle points were introduced and it was established that the
existence of an open-loop saddle point for the problem is equivalent to the solvability of a
forward-backward stochastic differential equation (FBSDE, for short), and the existence of a
closed-loop saddle point for the problem is equivalent to the solvability of a (differential) Riccati equation. In this paper, we will not assume \rf{16Apr2-1+2=0} so that \rf{16Apr2-zero} is not necessarily true. Such a Problem (SDG) is usually referred to as an LQ stochastic two-person {\it nonzero-sum} differential game, emphasizing that \rf{16Apr2-zero} is not assumed. We have two main goals in this paper: Establish a theory for Problem (SDG) parallel to that of \cite{Sun-Yong 2014} (for zero-sum case); and study the difference between the closed-loop representation of open-loop Nash equilibria and the outcome of closed-loop Nash equilibria. It turns out that the above-mentioned difference for the non-zero sum case is indicated through the symmetry of the corresponding Riccati equations: One is symmetric and the other is not.
On the other hand, we found that the situation in the zero-sum case, which was not discussed in \cite{Sun-Yong 2014}, is totally different: The closed-loop representation of open-loop saddle points coincides with the outcome of  the corresponding closed-loop saddle point, when both exist. In particular, for stochastic linear quadratic optimal control problem, the closed-loop representation of open-loop optimal controls is the outcome of the corresponding closed-loop optimal strategy (\cite{Sun-Li-Yong 2016}).

\ms

Mathematically, posing condition \rf{16Apr2-1+2=0} makes the structure of the problem much simpler, since with such a condition, only one performance index is needed, for which one player is the minimizer and the other player is the maximizer. However, as we know that in the real life, each player should have his/her own cost functional, and even for the totally hostile situation, the objectives of the opponents might not necessarily be exactly the opposite (zero-sum). Therefore, realistically, it is more meaningful to investigate Problem (SDG) without assuming \rf{16Apr2-1+2=0}. By the way, although we will not discuss such a situation in the current paper, we still would like to point out that sometimes, certain cooperations between the players might result in both players rewarded more.

\ms

Static version of nonzero-sum differential games could be regarded as a kind of non-cooperative games for which one can trace back to the work of Nash \cite{Nash 1951}. For some early works on nonzero-sum differential games, we would like to mention Lukes--Russell \cite{Lukes-Russell 1971}, Friedman \cite{Friedman 1972}, and Bensoussan \cite{Bensoussan 1974}. In the past two decays, due to the appearance of backward stochastic differential equations (BSDEs, for short), some new and interesting works published; Among them, we would like to mention \cite{Hamadene 1998, Hamadene 1999,El Karoui-Hamadene 2003, Buckdahn-Cardaliaguet-Rainer 2004, Rainer 2007, Hamadene-Mu 2015}.

\ms

The rest of the paper is organized as follows. Section 2 will collect some preliminaries.
Among other things, we will recall some known results on LQ optimal control problems.
In Section 3, we will introduce open-loop and closed-loop Nash equilibria. A characterization of the existence of open-loop Nash equilibria in terms of solvability of two coupled FBSDEs will be
presented in Section 4. Section 5 is devoted to the discussion on the closed-loop Nash equilibria whose existence is characterized by the solvability of two coupled symmetric Riccati equations. In Section 6, we will present two examples showing the difference between open-loop and closed-loop Nash equilibria. In Section 7, closed-loop representation of open-loop Nash equilibria will be studied, and comparison between
the closed-loop representation of open-loop Nash equilibria and the outcome of closed-loop Nash
equilibria will be carried out. Finally, we will take a deeper look at the situation for LQ zero-sum games in Section 8.

\section{Preliminaries}

Let $\dbR^{n\times m}$ be the space of all $(n\times m)$ matrices and $\dbS^n\subseteq\dbR^{n\times n}$ be
the set of all $(n\times n)$ symmetric matrices. The inner product $\lan\cd\,,\cd\ran$ on $\dbR^{n\times m}$
is given by $\lan M,N\ran\mapsto\tr(M^\top N)$, and the induced norm is given by $|M|=\sqrt{\tr(M^\top M)}$. We denote by $\sR(M)$ the range of
a matrix $M$, and for $M,N\in\dbS^n$ we use the notation $M\ges N$ (respectively, $M>N$) to indicate that
$M-N$ is positive semi-definite (respectively, positive definite).
Recall that any $M\in\dbR^{n\times m}$ admits a unique (Moore--Penrose) {\it pseudo-inverse}
$M^\dag\in\dbR^{m\times n}$ having the following properties (\cite{Penrose 1955}):
$$MM^\dag M=M, \q M^\dag MM^\dag=M^\dag, \q (MM^\dag)^\top=MM^\dag, \q (M^\dag M)^\top=M^\dag M.$$
Further, if $M\in\dbR^{n\times m}$ and $\Psi\in\dbR^{n\times\ell}$ such that
$$\sR(\Psi)\subseteq\sR(M),$$
then all the solutions $\Th$ to the linear equation
$$M\Th=\Psi$$
are given by the following:
$$\Th=M^\dag\Psi+(I-M^\dag M)\G,\qq\G\in\dbR^{m\times\ell}.$$
In addition, if $M=M^\top\in\dbS^n$, then
$$M^\dag=(M^\dag)^\top, \q MM^\dag=M^\dag M; \q\hb{and}\q M\ges0\iff M^\dag\ges0.$$
Next, let $T>0$ be a fixed time horizon. For any $t\in[0,T]$ and Euclidean space $\dbH$,
we introduce the following spaces of deterministic functions:
$$\ba{ll}
\ns\ds L^p(t,T;\dbH)=\Big\{\f:[t,T]\to\dbH\bigm|\int_t^T|\f(s)|^pds<\infty\Big\},\q 1\les p<\infty,\ea$$
$$\ba{ll}
\ns\ds L^\infty(t,T;\dbH)=\Big\{\f:[t,T]\to\dbH\bigm|\esssup_{s\in[t,T]}|\f(s)|<\infty\Big\},\\
\ns\ds C([t,T];\dbH)=\Big\{\f:[t,T]\to\dbH\bigm|\f(\cd)\hb{ is continuous}\Big\}.\ea$$
Further, we introduce the following spaces of random variables and stochastic processes: For any $t\in[0,T]$,
$$\ba{ll}
\ds L^2_{\cF_t}(\O;\dbH)=\Big\{\xi:\O\to\dbH\bigm|\xi\hb{ is $\cF_t$-measurable, }\dbE|\xi|^2<\i\Big\},\\
\ns\ds L_\dbF^2(t,T;\dbH)=\Big\{\f:[t,T]\times\O\to\dbH\bigm|\f(\cd)\hb{ is $\dbF$-progressively measurable, }
\dbE\int^T_t|\f(s)|^2ds<\infty\Big\},\\
\ns\ds L_\dbF^2(\O;C([t,T];\dbH))=\Big\{\f:[t,T]\times\O\to\dbH\bigm|\f(\cd)\hb{ is $\dbF$-adapted, continuous, }
\dbE\(\sup_{t\les s\les T}|\f(s)|^2\)<\infty\Big\},\\
\ns\ds L^2_\dbF(\O;L^1(t,T;\dbH))=\Big\{\f:[t,T]\times\O\to\dbH\bigm|\f(\cd)\hb{ is $\dbF$-progressively measurable, }
\dbE\(\int_t^T|\f(s)|ds\)^2<\infty\Big\}.\ea$$

\ss

We now recall some results on stochastic LQ optimal control problems. Consider the state equation
\bel{LQ-state}\left\{\2n\ba{ll}
\ds dX(s)=\big[A(s)X(s)+B(s)u(s)+b(s)\big]ds\\
\ns\ds\qq\qq~+\big[C(s)X(s)+D(s)u(s)+\si(s)\big]dW(s),\qq s\in[t,T],\\
\ns\ds X(t)=x. \ea\right.\ee
The cost functional takes the following form:
\bel{LQ-cost}\ba{ll}
\ds J(t,x;u(\cd))\deq\dbE\Big\{\lan GX(T),X(T)\ran+2\lan g,X(T)\ran\\
\ns\ds\qq\qq\qq\qq
+\int_t^T\[\lan{\sc\begin{pmatrix}\sc Q(s) & \sc S(s)^\top \\
                                  \sc S(s) & \sc R(s)      \end{pmatrix}
\begin{pmatrix}\sc X(s) \\ \sc u(s)\end{pmatrix},
\begin{pmatrix}\sc X(s) \\ \sc u(s)\end{pmatrix}}\ran
+2\lan{\sc\begin{pmatrix}\sc q(s) \\ \sc\rho(s)\end{pmatrix},
          \begin{pmatrix}\sc X(s) \\ \sc u(s)  \end{pmatrix}}\ran\]ds\Big\}. \ea\ee
We adopt the following assumptions.

\ms

{\bf(S1)} The coefficients of the state equation satisfy the following:
$$\left\{\2n\ba{llll}
A(\cd)\in L^1(0,T;\dbR^{n\times n}), &~ B(\cd)\in L^2(0,T;\dbR^{n\times m}),
&~ b(\cd)\in L^2_\dbF(\O;L^1(0,T;\dbR^n)),\\
\ns C(\cd)\in L^2(0,T;\dbR^{n\times n}), &~ D(\cd)\in L^\i(0,T;\dbR^{n\times m}),
&~ \si(\cd)\in L_\dbF^2(0,T;\dbR^n).
\ea\right.$$

{\bf(S2)} The weighting coefficients in the cost functional satisfy the following:
$$\left\{\2n\ba{llll}
Q(\cd)\in L^1(0,T;\dbS^n), & S(\cd)\in L^2(0,T;\dbR^{m\times n}),
& R(\cd)\in L^\i(0,T;\dbS^m),\\
\ns q(\cd)\in L^2_\dbF(\O;L^1(0,T;\dbR^n)), & \rho(\cd)\in L_\dbF^2(0,T;\dbR^m),
& g\in L^2_{\cF_T}(\O;\dbR^n), \q G\in\dbS^n.
\ea\right.$$

Note that under (S1), for any $(t,x)\in[0,T)\times\dbR^n$ and
$u(\cd)\in\cU[t,T]\equiv L^2_\dbF(t,T;\dbR^m)$, the state equation \rf{LQ-state} admits
a unique strong solution $X(\cd)\equiv X(\cd\,;t,x,u(\cd))$. Further, if (S2) is also assumed,
then the cost functional \rf{LQ-cost} is well-defined for every $(t,x)\in[0,T)\times\dbR^n$
and $u(\cd)\in\cU[t,T]$. Therefore, the following problem is meaningful.

\ms

\bf Problem (SLQ). \rm For any given initial pair $(t,x)\in[0,T)\times\dbR^n$,
find a $\bar u(\cd)\in\cU[t,T]$ such that
\bel{J(u)=inf}J(t,x;\bar u(\cd))=\inf_{u(\cd)\in\cU[t,T]}J(t,x;u(\cd)).\ee

Any $\bar u(\cd)\in\cU[t,T]$ satisfying \rf{J(u)=inf} is called an {\it open-loop optimal control} of Problem (SLQ) for $(t,x)$; the corresponding $\bar X(\cd)\equiv X(\cd\,;t,x,\bar u(\cd))$ is called an {\it open-loop optimal state process} and $(\bar X(\cd),\bar u(\cd))$ is called an {\it open-loop optimal pair}.

\bde{open-loop solvable} \rm Let $(t,x)\in[0,T)\times\dbR^n$. If there exists a (unique) $\bar u(\cd)\in\cU[t,T]$ such that \rf{J(u)=inf} holds,
then we say that Problem (SLQ) is ({\it uniquely}) {\it open-loop solvable at} $(t,x)$.
If Problem (SLQ) is (uniquely) open-loop solvable for any $(t,x)\in[0,T)\times\dbR^n$, then we say that Problem (SLQ) is
({\it uniquely}) {\it open-loop solvable on} $[0,T)\times\dbR^n$.
\ede

%
%

The following result is concerned with open-loop optimal controls of Problem (SLQ) for a given
initial pair, whose proof can be found in \cite{Sun-Yong 2014} (see also \cite{Sun-Li-Yong 2016}).

\bt{bt-16Apr2-17:10}\sl Let {\rm(S1)--(S2)} hold. For a given initial pair
$(t,x)\in[0,T)\times\dbR^n$, a state-control pair $(\bar X(\cd),\bar u(\cd))$
is an open-loop optimal pair of Problem {\rm(SLQ)} if and only if the following hold:

\ms

{\rm(i)} The {\it stationarity condition} holds:
$$B(s)^\top\bar Y(s)+D(s)^\top\bar Z(s)+S(s)\bar X(s)+R(s)\bar u(s)+\rho(s)=0,
\q \ae~s\in[t,T],~\as$$
where $(\bar Y(\cd),\bar Z(\cd))$ is the {\it adapted solution} to the following BSDE:
$$\left\{\2n\ba{ll}
\ds d\bar Y(s)=-\big[A(s)^\top\bar Y(s)+C(s)^\top\bar Z(s)+Q(s)\bar X(s)+S(s)^\top\bar u(s)+q(s)\big]ds+\bar Z(s)dW(s),\q s\in[t,T],\\
\ns\ds\bar Y(T)=G\bar X(T)+g. \ea\right.$$

{\rm(ii)} The map $u(\cd)\mapsto J(t,0;u(\cd))$ is convex.

\et

Next, for any given $t\in[0,T)$, take $\Th(\cd)\in L^2(t,T;\dbR^{m\times n})\equiv\cQ[t,T]$
and $v(\cd)\in\cU[t,T]$. For any $x\in\dbR^n$, let us consider the following equation:
\bel{closed-loop1}\left\{\2n\ba{ll}
\ds dX(s)=\big\{[A(s)+B(s)\Th(s)]X(s)+B(s)v(s)+b(s)\big\}ds\\
\ns\ds\qq\qq~+\big\{[C(s)+D(s)\Th(s)]X(s)+D(s)v(s)+\si(s)\big\}dW(s),\qq s\in[t,T],\\
\ns\ds X(t)=x, \ea\right.\ee
which admits a unique solution $X(\cd)\equiv X(\cd\,;t,x,\Th(\cd),v(\cd))$, depending on
$\Th(\cd)$ and $v(\cd)$. The above is called a {\it closed-loop system} of the original
state equation \rf{LQ-state} under {\it closed-loop strategy} $(\Th(\cd),v(\cd))$.
We point out that $(\Th(\cd),v(\cd))$ is independent of the initial state $x$.
With the above corresponding solution $X(\cd)$, we define
$$\ba{ll}
\ds J(t,x;\Th(\cd)X(\cd)+v(\cd))=\dbE\Big\{\lan GX(T),X(T)\ran+2\lan g,X(T)\ran\\
\ns\ds\qq\qq\qq~
+\int_t^T\[\lan{\sc\begin{pmatrix}\sc Q(s) &\2n \sc S(s)^\top \\
                                  \sc S(s) &\2n \sc R(s)      \end{pmatrix}
\begin{pmatrix}\sc X(s) \\ \sc\Th(s)X(s)+v(s) \end{pmatrix},
\begin{pmatrix}\sc X(s) \\ \sc\Th(s)X(s)+v(s)     \end{pmatrix}}\ran
+2\lan{\sc\begin{pmatrix}\sc q(s) \\ \sc\rho(s)    \end{pmatrix},
          \begin{pmatrix}\sc X(s) \\ \sc\Th(s)X(s)+v(s)\end{pmatrix}}\ran\]ds\Big\}. \ea$$
We now recall the following definition.

\bde{bde-16Apr2-21:00}\rm A pair $(\bar\Th(\cd),\bar v(\cd))\in\cQ[t,T]\times\cU[t,T]$
is called a {\it closed-loop optimal strategy} of Problem (SLQ) on $[t,T]$ if
\bel{16Apr2-21:30}\ba{ll}
\ds J(t,x;\bar\Th(\cd)\bar X(\cd)+\bar v(\cd))\les J(t,x;\Th(\cd)X(\cd)+v(\cd)),\\
\ns\ds\qq\qq\qq\qq\qq\qq\q\forall\, x\in\dbR^n,~\forall\,(\Th(\cd),v(\cd))\in\cQ[t,T]\times\cU[t,T], \ea\ee
where $\bar X(\cd)=X(\cd\,;t,x,\bar\Th(\cd),\bar v(\cd))$, and $X(\cd)=X(\cd\,;t,x,\Th(\cd),v(\cd))$.
\ede

We emphasize that the pair $(\bar\Th(\cd),\bar v(\cd))$ is required to be independent of the
initial state $x\in\dbR^n$. It is interesting that the following equivalent theorem holds.

\bp{strategy equivalence}\sl Let {\rm(S1)--(S2)} hold and let
$(\bar\Th(\cd),\bar v(\cd))\in\cQ[t,T]\times\cU[t,T]$. Then the following statements are equivalent:

\ms

{\rm(i)} $(\bar\Th(\cd),\bar v(\cd))$ is a closed-loop optimal strategy of Problem {\rm(SLQ)} on $[t,T]$.

\ms

{\rm(ii)} For any $x\in\dbR^n$ and $v(\cd)\in\cU[t,T]$,
$$J(t,x;\bar\Th(\cd)\bar X(\cd)+\bar v(\cd))\les J(t,x;\bar\Th(\cd)X(\cd)+v(\cd)),$$
where $\bar X(\cd)=X(\cd\,;t,x,\bar\Th(\cd),\bar v(\cd))$ and $X(\cd)=X(\cd\,;t,x,\bar\Th(\cd),v(\cd))$.

\ms

{\rm(iii)} For any $x\in\dbR^n$ and $u(\cd)\in\cU[t,T]$,
\bel{bar v<u}J(t,x;\bar\Th(\cd)\bar X(\cd)+\bar v(\cd))\les J(t,x;u(\cd)),\ee
where $\bar X(\cd)=X(\cd\,;t,x,\bar\Th(\cd),\bar v(\cd))$.

\ep

\it Proof. \rm The implication (i) $\Ra$ (ii) follows by taking $\Th(\cd)=\bar\Th(\cd)$
in \rf{16Apr2-21:30}.

\ms

For the implication (ii) $\Ra$ (iii), take any $u(\cd)\in\cU[t,T]$ and let
$X(\cd)=X(\cd\,;t,x,u(\cd))$. Then
$$\ba{ll}
\ds dX(s)=\big\{[A(s)+B(s)\bar\Th(s)]X(s)+B(s)[u(s)-\bar\Th(s)X(s)]+b(s)\big\}ds\\
\ns\ds\qq\qq\q
+\big\{[C(s)+D(s)\bar\Th(s)]X(s)+D(s)[u(s)-\bar\Th(s)X(s)]+\si(s)\big\}dW(s),\ea$$
with $X(t)=x$. Thus, if let
$$v(\cd)=u(\cd)-\bar\Th(\cd)X(\cd),$$
we have
$$J(t,x;\bar\Th(\cd)\bar X(\cd)+\bar v(\cd))\les J(t,x;\bar\Th(\cd)X(\cd)+v(\cd))=J(t,x;u(\cd)),$$
which proves (iii).

\ms

For the implication (iii) $\Ra$ (i), take any $(\Th(\cd),v(\cd))\in\cQ[t,T]\times\cU[t,T]$
and let $X(\cd)$ be the solution to \rf{closed-loop1}. Let $u(\cd)=\Th(\cd)X(\cd)+v(\cd)$, Then by (iii), we have
$$J(t,x;\bar\Th(\cd)\bar X(\cd)+\bar v(\cd))\les J(t,x;u(\cd))=J(t,x;\Th(\cd)X(\cd)+v(\cd)).$$
This completes the proof. \endpf

\ms

From the above result, we see that if $(\bar\Th(\cd),\bar v(\cd))$ is a closed-loop optimal
strategy of Problem (SLQ) on $[t,T]$, then for any fixed initial state $x\in\dbR^n$,
with $\bar X(\cd)$ denoting the state process corresponding to $(t,x)$ and $(\bar\Th(\cd),\bar v(\cd))$,
\rf{bar v<u} implies that the outcome
$$\bar u(\cd)\equiv\bar\Th(\cd)\bar X(\cd)+\bar v(\cd)\in\cU[t,T]$$
is an open-loop optimal control of Problem (SLQ) for $(t,x)$. Therefore, for Problem (SLQ),
the existence of closed-loop strategies on $[t,T]$ implies the existence of open-loop optimal controls for initial pair $(t,x)$ for any $x\in\dbR^n$. We point out that the situation will be different for two-person differential games. Details will be carried out later.

\ms

For closed-loop optimal strategies, we have the following characterization
(\cite{Sun-Yong 2014,Sun-Li-Yong 2016}).

\bt{closed-loop strategy}\sl Let {\rm(S1)--(S2)} hold. Then Problem {\rm(SLQ)} admits a
closed-loop optimal strategy on $[t,T]$ if and only if the following Riccati equation
admits a solution $P(\cd)\in C([t,T];\dbS^n)$:
$$\left\{\2n\ba{ll}
\ds\dot P+PA+A^\top P+C^\top PC+Q\\
\ns\ds\,~~-(PB+C^\top PD+S^\top)(R+D^\top PD)^\dag(B^\top P+D^\top PC+S)=0,\qq\ae~\hb{on }[t,T],\\
\ns\ds\sR(B^\top P+D^\top PC+S)\subseteq\sR(R+D^\top PD),\qq\ae~\hb{on }[t,T],\\
\ns\ds R+D^\top PD\ges0,\qq\ae~\hb{on }[t,T],\\
\ns\ds P(T)=G, \ea\right.$$
such that
$$(R+D^\top PD)^\dag(B^\top P+D^\top PC+S)\in L^2(t,T;\dbR^{m\times n}),$$
and the adapted solution $(\eta(\cd),\z(\cd))$ to the BSDE
$$\left\{\2n\ba{ll}
\ds d\eta=-\,\Big\{\big[A-B(R+D^\top PD)^\dag(B^\top P+D^\top PC+S)\big]^\top\eta\\
\ns\ds\qq\q~~+\big[C-D(R+D^\top PD)^\dag(B^\top P+D^\top PC+S)\big]^\top\z\\
\ns\ds\qq\q~~+\big[C-D(R+D^\top PD)^\dag(B^\top P+D^\top PC+S)\big]^\top P\si\\
\ns\ds\qq\q~~-(PB+C^\top PD+S^\top)(R+D^\top PD)^\dag\rho+Pb+q\Big\}ds+\z dW,\q~ s\in[t,T],\\
\ns\ds\eta(T)=g, \ea\right.$$
satisfies
$$\left\{\2n\ba{ll}
\ds B^\top\eta+D^\top\z+D^\top P\si+\rho\in\sR(R+D^\top PD),\q \ae~s\in[t,T],~\as\\
\ns\ds (R+D^\top PD)^\dag(B^\top\eta+D^\top\z+D^\top P\si+\rho)\in L_\dbF^2(t,T;\dbR^m).
\ea\right.$$
In this case, any closed-loop optimal strategy $(\bar\Th(\cd),\bar v(\cd))$
of Problem {\rm(SLQ)} admits the following representation:
$$\left\{\2n\ba{cll}
\bar\Th \3n&=&\2n\5n~ -(R+D^\top PD)^\dag(B^\top P+D^\top PC+S)
+\big[I-(R+D^\top PD)^\dag(R+D^\top PD)\big]\th,\\
\ns\bar v \3n&=&\2n\5n~ -(R+D^\top PD)^\dag(B^\top\eta+D^\top\z+D^\top P\si+\rho)
+\big[I-(R+D^\top PD)^\dag(R+D^\top PD)\big]\n,
\ea\right.$$
for some $\th(\cd)\in L^2(t,T;\dbR^{m\times n})$ and $\n(\cd)\in L_\dbF^2(t,T;\dbR^m)$.
Further, the value function is given by
$$\ba{ll}
\ns\ds V(t,x)=\dbE\Big\{\lan P(t)x,x\ran+2\lan\eta(t),x\ran
+\int_t^T\[\lan P\si,\si\ran+2\lan\eta,b\ran+2\lan\z,\si\ran\\
\ns\ds\qq\qq\q~-\blan(R+D^\top PD)^\dag(B^\top\eta+D^\top\z+D^\top P\si+\rho),
B^\top\eta+D^\top\z+D^\top P\si+\rho\bran\]ds\Big\}. \ea$$
\et

\section{Stochastic Differential Games}

We return to our Problem (SDG). Recall the sets $\cU_i[t,T]=L^2_\dbF(t,T;\dbR^{m_i})$
of all open-loop controls of Player $i$ ($i=1,2$).
For notational simplicity, we let $m=m_1+m_2$ and denote
$$\ba{ll}
\ds B(\cd)=(B_1(\cd),B_2(\cd)),\q D(\cd)=(D_1(\cd),D_2(\cd)),\\
\ns\ds S^i(\cd)\1n=\1n\begin{pmatrix}S^i_1(\cd) \\ S^i_2(\cd)\end{pmatrix},
    \q R^i(\cd)\1n=\1n\begin{pmatrix}R^i_{11}(\cd) & R^i_{12}(\cd)\\
                                     R^i_{21}(\cd) & R^i_{22}(\cd)\end{pmatrix}
          \1n\equiv\1n\begin{pmatrix}R^i_1(\cd) \\ R^i_2(\cd)\end{pmatrix},
    \q \rho^i(\cd)\1n=\1n\begin{pmatrix}\rho^i_1(\cd) \\ \rho^i_2(\cd)\end{pmatrix},
    \q u(\cd)\1n=\1n\begin{pmatrix}u_1(\cd) \\ u_2(\cd)\end{pmatrix}. \ea$$
Naturally, we identify $\cU[t,T]=\cU_1[t,T]\times\cU_2[t,T]$.
With such notations, the state equation becomes
\bel{16Apr3-state}\left\{\2n\ba{ll}
\ds dX(s)=\big[A(s)X(s)+B(s)u(s)+b(s)\big]ds\\
\ns\ds\qq\qq\2n~+\big[C(s)X(s)+D(s)u(s)+\si(s)\big]dW(s), \qq s\in[t,T],\\
\ns\ds X(t)=x, \ea\right.\ee
and the cost functionals become ($i=1,2$)
%
$$\ba{ll}
\ds J^i(t,x;u(\cd))=\dbE\Big\{\lan G^iX(T),X(T)\ran+2\lan g^i,X(T)\ran\\
\ns\ds\qq\qq\qq\qq+\int_t^T\[\blan
\sc{\begin{pmatrix}\sc Q^i(s) &\sc S^i(s)^\top \\ \sc S^i(s) &\sc R^i(s)\end{pmatrix}
\begin{pmatrix}\sc X(s)   \\ \sc u(s)      \end{pmatrix},
\begin{pmatrix}\sc X(s)   \\ \sc u(s)      \end{pmatrix}}\bran+2\blan
\begin{pmatrix}\sc q^i(s) \\ \sc\rho^i(s) \end{pmatrix},
\begin{pmatrix}\sc X(s)   \\ \sc u(s)      \end{pmatrix}\bran\]ds\Big\}. \ea$$
Now let us introduce the following standard assumptions:

\ms

{\bf(G1)} The coefficients of the state equation satisfy the following:
$$\left\{\2n\ba{lll}
\ds A(\cd)\in L^1(0,T;\dbR^{n\times n}), & B(\cd)\in L^2(0,T;\dbR^{n\times m}), & b(\cd)\in L^2_\dbF(\O;L^1(0,T;\dbR^n)),\\
\ns\ds C(\cd)\in L^2(0,T;\dbR^{n\times n}), & D(\cd)\in L^\i(0,T;\dbR^{n\times m}), & \si(\cd)\in L_\dbF^2(0,T;\dbR^n).
\ea\right.$$

{\bf(G2)} The weighting coefficients in the cost functionals satisfy the following: For $i=1,2$,
$$\left\{\2n\ba{ll}
\ds Q^i(\cd)\in L^1(0,T;\dbS^n),\q S^i(\cd)\in L^2(0,T;\dbR^{m\times n}),
\q R^i(\cd)\in L^\i(0,T;\dbS^m),\\
\ns\ds q^i(\cd)\in L^2_\dbF(\O;L^1(0,T;\dbR^n)),\q\rho^i(\cd)\in L_\dbF^2(0,T;\dbR^m),\q
g^i\in L^2_{\cF_T}(\O;\dbR^n),\q G^i\in\dbS^n.\ea\right.$$

Under (G1), for any $(t,x)\in[0,T)\times\dbR^n$ and
$u(\cd)=(u_1(\cd)^\top,u_2(\cd)^\top)^\top\in\cU[t,T]$, equation \rf{16Apr3-state}
admits a unique solution (\cite{Yong-Zhou 1999})
$$X(\cd)\deq X(\cd\,;t,x,u_1(\cd),u_2(\cd))\equiv X(\cd\,;t,x,u(\cd))
\in L^2_\dbF(\O;C([t,T];\dbR^n)).$$
Moreover, the following estimate holds:
$$\dbE\(\sup_{t\les s\les T}|X(s)|^2\)\les
K\dbE\Big\{|x|^2+\(\int_t^T|b(s)|ds\)^2+\int_t^T|\si(s)|^2ds+\int^T_t|u(s)|^2ds\Big\},$$
where $K>0$ represents a generic constant. Therefore, under (G1)--(G2), the cost functionals
$J^i(t,x;u(\cd))\equiv J^i(t,x;u_1(\cd),u_2(\cd))$ are well-defined for all
$(t,x)\in[0,T)\times\dbR^n$ and all $(u_1(\cd),u_2(\cd))\in\cU_1[t,T]\times\cU_2[t,T]$.
Having the above, we now introduce the following definition.

\bde{open-Nash-equilibrium} \rm A pair $(u^*_1(\cd),u^*_2(\cd))\in\cU_1[t,T]\times\cU_2[t,T]$
is called an {\it open-loop Nash equilibrium} of Problem (SDG) for the initial pair
$(t,x)\in[0,T)\times\dbR^n$ if
\bel{Nash-open*}\ba{ll}
\ns\ds J^1(t,x;u^*_1(\cd),u^*_2(\cd))\les J^1(t,x;u_1(\cd),u^*_2(\cd)),\qq\forall\, u_1(\cd)\in\cU_1[t,T],\\
\ns\ds J^2(t,x;u^*_1(\cd),u^*_2(\cd))\les J^2(t,x;u^*_1(\cd),u_2(\cd)),
\qq\forall\, u_2(\cd)\in\cU_2[t,T].\ea\ee

\ede

Next, we denote
$$\cQ_i[t,T]=L^2(t,T;\dbR^{m_i\times n}),\qq i=1,2.$$
For any initial pair $(t,x)\in[0,T)\times\dbR^n$,
$\Th(\cd)\equiv(\Th_1(\cd)^\top,\Th_2(\cd)^\top)^\top\in\cQ_1[t,T]\times\cQ_2[t,T]$
and any $v(\cd)\equiv(v_1(\cd)^\top$, $v_2(\cd)^\top)^\top\in\cU_1[t,T]\times\cU_2[t,T]$,
consider the following system:
\bel{16Apr4-state-closed}\left\{\2n\ba{ll}
\ds dX(s)=\big\{[A(s)+B(s)\Th(s)]X(s)+B(s)v(s)+b(s)\big\}ds\\
\ns\ds\qq\qq\1n~+\big\{[C(s)+D(s)\Th(s)]X(s)+D(s)v(s)+\si(s)\big\}dW(s),\qq s\in[t,T],\\
\ns\ds X(t)= x. \ea\right.\ee
Under (G1), the above admits a unique solution
$X(\cd)\equiv X(\cd\,;t,x,\Th_1(\cd),v_1(\cd),\Th_2(\cd),v_2(\cd))$. If we denote
\bel{ui}u_i(\cd)=\Th_i(\cd)X(\cd)+v_i(\cd),\qq i=1,2,\ee
then the above \rf{16Apr4-state-closed} coincides with the original state equation \rf{state}.
We call $(\Th_i(\cd),v_i(\cd))$ a {\it closed-loop strategy} of Player $i$, and call
\rf{16Apr4-state-closed} the {\it closed-loop system} of the original system under closed-loop
strategies $(\Th_1(\cd),v_1(\cd))$ and $(\Th_2(\cd),v_2(\cd))$ of Players 1 and 2. Also, we call
$u(\cd)\equiv(u_1(\cd)^\top,u_2(\cd)^\top)^\top$ with $u_i(\cd)$ defined by \rf{ui} the outcome of
the closed-loop strategy $(\Th(\cd),v(\cd))$. With the solution $X(\cd)$ to \rf{16Apr4-state-closed},
we denote
\bel{16Apr7-22:00}\ba{ll}
\no\ms\ds J^i(t,x;\Th(\cd)X(\cd)+v(\cd))\equiv
J^i(t,x;\Th_1(\cd)X(\cd)+v_1(\cd),\Th_2(\cd)X(\cd)+v_2(\cd))\\
\no\ms\ds=\dbE\Big\{\lan G^iX(T),X(T)\ran+2\lan g^i,X(T)\ran\\
\no\ms\ds\q~+\int_t^T\[\blan{\sc
\begin{pmatrix}\sc Q^i & \sc (S^i)^\top\\
               \sc S^i & \sc R^i       \end{pmatrix}
\begin{pmatrix}\sc X   \\ \sc\Th X+v   \end{pmatrix},
\begin{pmatrix}\sc X   \\ \sc\Th X+v   \end{pmatrix}}\bran
+2\blan{\sc
\begin{pmatrix}\sc q^i \\ \sc\rho^i \end{pmatrix},
\begin{pmatrix}\sc X   \\ \sc\Th X+v\end{pmatrix}}\bran\]ds\Big\}\\
\no\ms\ds=\dbE\Big\{\lan G^iX(T),X(T)\ran+2\lan g^i,X(T)\ran\\
\ds\q~+\int_t^T\[\blan{\sc
\begin{pmatrix}\sc Q^i+\Th^\top S^i+(S^i)^\top\Th+\Th^\top R^i\Th &\sc (S^i)^\top+\Th^\top R^i\\
               \sc S^i+R^i\Th &\sc R^i \end{pmatrix}
\begin{pmatrix}\sc X \\ \sc v\end{pmatrix},
\begin{pmatrix}\sc X \\ \sc v\end{pmatrix}}\bran
+2\blan{\sc
\begin{pmatrix}\sc q^i+\Th^\top\rho^i \\ \sc\rho^i\end{pmatrix},
\begin{pmatrix}\sc X \\ \sc v \end{pmatrix}}\bran\]ds\Big\}.\ea\ee
Similarly, one can define $J^i(t,x;\Th_1(\cd)X(\cd)+v_1(\cd),u_2(\cd))$ and
$J^i(t,x;u_1(\cd),\Th_2(\cd)X(\cd)+v_2(\cd))$. We now introduce the following definition.

\bde{bde-16Apr4-17:00}\rm A 4-tuple $(\Th_1^*(\cd),v_1^*(\cd);\Th_2^*(\cd),v_2^*(\cd))
\in\cQ_1[t,T]\times\cU_1[t,T]\times\cQ_2[t,T]\times\cU_2[t,T]$ is called a
{\it closed-loop Nash equilibrium} of Problem (SDG) on $[t,T]$ if for any $x\in\dbR^n$ and
any 4-tuple $(\Th_1(\cd),v_1(\cd);$ $\Th_2(\cd),v_2(\cd))\in
\cQ_1[t,T]\times\cU_1[t,T]\times\cQ_2[t,T]\times\cU_2[t,T]$, the following hold:
\begin{eqnarray}
&&\label{Nash-closed1}\ba{lll}
\ds J^1(t,x;\Th_1^*(\cd)X^*(\cd)+v_1^*(\cd),\Th_2^*(\cd)X^*(\cd)+v_2^*(\cd))\les J^1(t,x;\Th_1(\cd)X(\cd)+v_1(\cd),\Th_2^*(\cd)X(\cd)+v_2^*(\cd)),\ea\\
\ns&&\label{Nash-closed2}\ba{ll}
\ds J^2(t,x;\Th_1^*(\cd)X^*(\cd)+v_1^*(\cd),\Th_2^*(\cd)X^*(\cd)+v_2^*(\cd))\les J^2(t,x;\Th^*_1(\cd)X(\cd)+v^*_1(\cd),\Th_2(\cd)X(\cd)+v_2(\cd)).\ea
\end{eqnarray}
\ede

Note that in both \rf{Nash-closed1} and \rf{Nash-closed2},
$$X^*(\cd)=X(\cd\,;t,x,\Th_1^*(\cd),v_1^*(\cd),\Th_2^*(\cd),v_2^*(\cd)),$$
whereas, in \rf{Nash-closed1},
$$X(\cd)=X(\cd\,;t,x,\Th_1(\cd),v_1(\cd),\Th_2^*(\cd),v_2^*(\cd)),$$
and in \rf{Nash-closed2},
$$X(\cd)=X(\cd\,;t,x,\Th_1^*(\cd),v_1^*(\cd),\Th_2(\cd),v_2(\cd)).$$
Thus, $X(\cd)$ appeared in \rf{Nash-closed1} and \rf{Nash-closed2} are different in general.
We emphasize that the closed-loop Nash equilibrium $(\Th^*_1(\cd),v^*_1(\cd);\Th^*_2(\cd),v^*_2(\cd))$
is independent of the initial state $x$.
The following result provides some equivalent definitions of closed-loop Nash equilibrium.

\bp{bp-16Apr4-17:30}\sl Let {\rm(G1)--(G2)} hold and let $(\Th^*_1(\cd),v_1^*(\cd);
\Th^*_2(\cd),v^*_2(\cd))\in\cQ_1[t,T]\times\cU_1[t,T]\times\cQ_2[t,T]\times\cU_2[t,T]$.
Then the following are equivalent:

\ms

{\rm(i)} $(\Th^*_1(\cd),v_1^*(\cd);\Th^*_2(\cd),v^*_2(\cd))$ is a closed-loop Nash equilibrium
      of Problem {\rm(SDG)} on $[t,T]$.

\ms

{\rm(ii)} For any $(v_1(\cd),v_2(\cd))\in\cU_1[t,T]\times\cU_2[t,T]$,
$$\ba{ll}
\no\ss\ds J^1(t,x;\Th_1^*(\cd)X^*(\cd)+v_1^*(\cd),\Th_2^*(\cd)X^*(\cd)+v_2^*(\cd))\les J^1(t,x;\Th_1^*(\cd)X(\cd)+v_1(\cd),\Th_2^*(\cd)X(\cd)+v_2^*(\cd)),\\
\no\ss\ds J^2(t,x;\Th_1^*(\cd)X^*(\cd)+v_1^*(\cd),\Th_2^*(\cd)X^*(\cd)+v_2^*(\cd))\les J^2(t,x;\Th^*_1(\cd)X(\cd)+v^*_1(\cd),\Th_2^*(\cd)X(\cd)+v_2(\cd)).\ea$$

{\rm(iii)} For any $(u_1(\cd),u_2(\cd))\in\cU_1[t,T]\times\cU_2[t,T]$,
\begin{eqnarray}
&&\label{16Apr4-Nash-closed1*}\ba{ll}
\ds J^1(t,x;\Th_1^*(\cd)X^*(\cd)+v_1^*(\cd),\Th_2^*(\cd)X^*(\cd)+v_2^*(\cd))\les J^1(t,x;u_1(\cd),\Th_2^*(\cd)X(\cd)+v_2^*(\cd)),\ea\\
\ns&&\label{16Apr4-Nash-closed2*}\ba{ll}
\ds J^2(t,x;\Th_1^*(\cd)X^*(\cd)+v_1^*(\cd),\Th_2^*(\cd)X^*(\cd)+v_2^*(\cd))\les J^2(t,x;\Th^*_1(\cd)X(\cd)+v^*_1(\cd),u_2(\cd)).\ea
\end{eqnarray}
\ep

\it Proof. \rm The proof is similar to that of Proposition \ref{strategy equivalence}. \endpf

\ms

If we denote
\bel{16Apr4-bar u}\bar u_i(\cd)=\Th^*_i(\cd)X^*(\cd)+v^*_i(\cd),\qq i=1,2,\ee
then \rf{16Apr4-Nash-closed1*}--\rf{16Apr4-Nash-closed2*} become
\begin{eqnarray}
&&\label{16Apr4-20:00-1}J^1(t,x;\bar u_1(\cd),\bar u_2(\cd))\les
J^1(t,x;u_1(\cd),\Th_2^*(\cd)X(\cd)+v_2^*(\cd)),\\
\ns&&\label{16Apr4-20:00-2}J^2(t,x;\bar u_1(\cd),\bar u_2(\cd))\les
J^2(t,x;\Th^*_1(\cd)X(\cd)+v^*_1(\cd),u_2(\cd)).
\end{eqnarray}
Since in \rf{16Apr4-20:00-1}, $X(\cd)$ corresponds to $u_1(\cd)$ and
$(\Th_2^*(\cd),v_2^*(\cd))$, one might not have
$$\bar u_2(\cd)=\Th^*_2(\cd)X(\cd)+v_2^*(\cd).$$
Likewise, one might not have the following either:
$$\bar u_1(\cd)=\Th_1^*(\cd)X(\cd)+v_1^*(\cd).$$
Hence, comparing this with \rf{Nash-open*}, we see that the outcome $(\bar u_1(\cd),\bar u_2(\cd))$
of the closed-loop Nash equilibrium $(\Th^*_1(\cd),v_1^*(\cd);\Th^*_2(\cd),v_2^*(\cd))$ defined
by \rf{16Apr4-bar u} is not an open-loop Nash equilibrium of Problem (SDG) for $(t,X^*(t))$ in general.

\ms

On the other hand, if $(\Th_1^*(\cd),v_1^*(\cd);\Th^*_2(\cd),v_2^*(\cd))$ is a closed-loop
Nash equilibrium of Problem (SDG) on $[t,T]$, we may consider the following state equation
(denoting $\Th^*(\cd)=(\Th_1^*(\cd)^\top\1n,\Th_2^*(\cd)^\top)^\top$)
\bel{16Apr4-20:30}\left\{\2n\ba{ll}
\ds dX(s)\1n=\1n\big[(A\1n+\1n B\Th^*)X\1n+\1n B_1v_1\1n+\1n B_2v_2\1n+\1n b\big]ds\1n+\1n\big[(C\1n+\1n D\Th^*)X\1n+\1n D_1v_1\1n+\1n D_2v_2\1n+\1n\si\big]dW(s),\q s\in[t,T],\\
\ns\ds X(t)=x, \ea\right.\ee
with cost functionals
\bel{16Apr4-20:40} \wt J^i(t,x;v_1(\cd),v_2(\cd))
=J^i(t,x;\Th^*_1(\cd)X(\cd)+v_1(\cd),\Th^*_2(\cd)X(\cd)+v_2(\cd)),\q~ i=1,2.\ee
Then by (ii) of Proposition \ref{bp-16Apr4-17:30}, $(v_1^*(\cd),v^*_2(\cd))$ is an open-loop
Nash equilibrium of the corresponding (nonzero-sum differential) problem.
Such an observation will be very useful below.

\section{Open-Loop Nash Equilibria and FBSDEs}

In this section, we discuss the open-loop Nash equilibria for Problem (SDG) in terms of FBSDEs.
The main result of this section can be stated as follows.

\bt{Theorem 4.1}\sl Let {\rm(G1)--(G2)} hold and let $(t,x)\in[0,T)\times\dbR^n$ be given.
Then $u^*(\cd)\equiv(u_1^*(\cd)^\top\1n,u_2^*(\cd)^\top)^\top\in\cU_1[t,T]\times\cU_2[t,T]$
is an open-loop Nash equilibrium of Problem {\rm(SDG)} for $(t,x)$ if and only if the following
two conditions hold:

\ms

{\rm(i)} For $i=1,2$, the adapted solution $(X^*(\cd),Y_i^*(\cd),Z_i^*(\cd))$ to the FBSDE on $[t,T]$
\bel{FBSDEi}\left\{\2n\ba{ll}
\ds dX^*(s)=\big[A(s)X^*(s)+B(s)u^*(s)+b(s)\big]ds+\big[C(s)X^*(s)+D(s)u^*(s)+\si(s)\big]dW(s),\\
\ns\ds dY_i^*(s)=-\big[A(s)^\top Y_i^*(s)+C(s)^\top Z_i^*(s)+Q^i(s)X^*(s)+S^i(s)^\top u^*(s)+q^i(s)\big]ds+Z_i^*(s)dW(s),\\
\ns\ds X^*(t)=x,\qq Y_i^*(T)=G^iX^*(T)+g^i, \ea\right.\ee
satisfies the following stationarity condition:
\bel{stationary}\ba{rl}
\ds B_i(s)^\top Y_i^*(s)+D_i(s)^\top Z_i^*(s)+S_i^i(s)X^*(s)+R^i_i(s)u^*(s)+\rho_i^i(s)=0,\q\ae~s\in[t,T],~\as\ea\ee

{\rm(ii)} For $i=1,2$, the following convexity condition holds:
\bel{convexity}\ba{ll}
\ds\dbE\Big\{\int_t^T\[\blan Q^i(s)X_i(s),X_i(s)\bran+2\blan S^i_i(s)X_i(s),u_i(s)\bran
+\blan R^i_{ii}(s)u_i(s),u_i(s)\bran\]ds\\
\ns\ds\qq\qq\q+\,\blan G^iX_i(T),X_i(T)\bran\Big\}\ges0,\qq\forall\,u_i(\cd)\in\cU_i[t,T],\ea\ee
where $X_i(\cd)$ is the solution to the following FSDE:
\bel{homogeneous}\left\{\2n\ba{ll}
\ds dX_i(s)=\big[A(s)X_i(s)+B_i(s)u_i(s)\big]ds+\big[C(s)X_i(s)+D_i(s)u_i(s)\big]dW(s),\qq s\in[t,T],\\
\ns\ds X_i(t)=0.\ea\right.\ee
Or, equivalently, the map $u_i(\cd)\mapsto J^i(t,x;u(\cd))$ is convex (for $i=1,2$).

\et

\it Proof. \rm For a given $(t,x)\in[0,T)\times\dbR^n$ and $u^*(\cd)\in\cU[t,T]$,
let $(X^*(\cd),Y_1^*(\cd),Z_1^*(\cd))$ be the adapted solution to FBSDE \rf{FBSDEi} with $i=1$.
For any $u_1(\cd)\in\cU_1[t,T]$ and $\e\in\dbR$, let $X^\e(\cd)$ be the solution to the following
perturbed state equation on $[t,T]$:
$$\left\{\2n\ba{ll}
\ds dX^\e(s)=\big\{A(s)X^\e(s)+B_1(s)[u_1^*(s)+\e u_1(s)]+B_2(s)u_2^*(s)+b(s)\big\}ds\\
\ns\ds\qq\qq~+\big\{C(s)X^\e(s)+D_1(s)[u_1^*(s)+\e u_1(s)]+D_2(s)u_2^*(s)+\si(s)\big\}dW(s),\\
\ns\ds X^\e(t)=x.\ea\right.$$
Then denoting $X_1(\cd)$ the solution of \rf{homogeneous} with $i=1$, we have $X^\e(\cd)=X^*(\cd)+\e X_1(\cd)$ and
$$\ba{ll}
\ds J^1(t,x;u^*_1(\cd)+\e u_1(\cd),u_2^*(\cd))-J^1(t,x;u_1^*(\cd),u_2^*(\cd))\\
\ns\ds=\e\dbE\Big\{\blan G^1[2X^*(T)+\e X_1(T)],X_1(T)\ran+2\lan g^1,X_1(T)\bran\\
\ns\ds\qq\q~+\int_t^T\[\blan{\sc
\begin{pmatrix}\sc Q^1   &\sc (S^1_1)^\top &\sc (S^1_2)^\top\\
               \sc S^1_1 &\sc R^1_{11}     &\sc R^1_{12}    \\
               \sc S^1_2 &\sc R^1_{21}     &\sc R^1_{22}    \end{pmatrix}
\begin{pmatrix}\sc 2X^*+\e X_1 \\ \sc 2u_1^*+\e u_1 \\ \sc 2u_2^* \end{pmatrix},
\begin{pmatrix}\sc X_1 \\ \sc u_1 \\ \sc 0 \end{pmatrix}}\bran
+2\blan{\sc\begin{pmatrix}\sc q^1 \\ \sc\rho^1_1\end{pmatrix},
           \begin{pmatrix}\sc X_1 \\ \sc u_1    \end{pmatrix}}\bran\]ds\Big\}\\
\ns\ds=2\e\dbE\Big\{\blan G^1X^*(T)+g^1,X_1(T)\bran+\int_t^T\[\blan Q^1X^*+(S^1)^\top u^*+q^1,X_1\bran
+\blan S^1_1X^*+R^1_1u^*+\rho^1_1,u_1\bran\]ds\Big\}\\
\ns\ds~~+\e^2\dbE\Big\{\blan G^1X_1(T),X_1(T)\bran+\int_t^T\[\blan Q^1X_1,X_1\bran
+2\blan S^1_1X_1,u_1\bran+\blan R^1_{11}u_1,u_1\bran\]ds\Big\}. \ea$$
On the other hand, applying It\^o's formula to $s\mapsto\lan Y_1^*(s),X_1(s)\ran$, we obtain
$$\ba{ll}
\ds\dbE\Big\{\blan G^1X^*(T)+g^1,X_1(T)\bran
+\int_t^T\[\blan Q^1X^*+(S^1)^\top u^*+q^1,X_1\bran
+\blan S^1_1X^*+R^1_1u^*+\rho^1_1,u_1\bran\]ds\Big\}\\
\ns\ds=\dbE\int_t^T\Big\{\blan-\big[A^\top Y_1^*+C^\top Z_1^*+Q^1X^*+(S^1)^\top u^*+q^1\big],X_1\bran
+\blan Y_1^*,AX_1+B_1u_1\bran\\
\ns\ds\qq\qq~+\blan Z_1^*,CX_1+D_1u_1\bran+\blan Q^1X^*+(S^1)^\top u^*+q^1,X_1\bran
+\blan S^1_1X^*+R^1_1u^*+\rho^1_1,u_1\bran\Big\}ds\\
\ns\ds=\dbE\int_t^T\blan B_1^\top Y_1^*+D_1^\top Z_1^*+S^1_1X^*+R^1_1u^*+\rho^1_1,u_1\bran ds.
\ea$$
Hence,
$$\ba{ll}
\ds J^1(t,x;u^*_1(\cd)+\e u_1(\cd),u_2^*(\cd))-J^1(t,x;u_1^*(\cd),u_2^*(\cd))\\
\ns\ds=2\e\dbE\int_t^T\blan B_1^\top Y_1^*+D_1^\top Z_1^*+S^1_1X^*+R^1_1u^*+\rho^1_1,u_1\bran ds\\
\ns\ds~~+\e^2\dbE\Big\{\blan G^1X_1(T),X_1(T)\bran+\int_t^T\[\blan Q^1X_1,X_1\bran
+2\blan S^1_1X_1,u_1\bran+\blan R^1_{11}u_1,u_1\bran\]ds\Big\}.\ea$$
It follows that
$$J^1(t,x;u_1^*(\cd),u_2^*(\cd))\les J^1(t,x;u_1^*(\cd)+\e u_1(\cd),u_2^*(\cd)),
\qq\forall\,u_1(\cd)\in\cU_1[t,T],~\forall\,\e\in\dbR,$$
if and only if \rf{convexity} holds for $i=1$, and
\bel{stationary1}B_1^\top Y_1^*+D_1^\top Z_1^*+S^1_1X^*+R^1_1u^*+\rho^1_1=0,
\qq\ae~s\in[t,T],~\as\ee
Similarly,
$$J^2(t,x;u_1^*(\cd),u_2^*(\cd))\les J^2(t,x;u_1^*(\cd),u_2^*(\cd)+\e u_2(\cd)),
\qq\forall\,u_2(\cd)\in\cU_2[t,T],~\forall\,\e\in\dbR,$$
if and only if \rf{convexity} holds for $i=2$, and
\bel{stationary2}B_2^\top Y_2^*+D_2^\top Z_2^*+S^2_2X^*+R^2_2u^*+\rho^2_2=0,
\qq\ae~s\in[t,T],~\as\ee
Combining \rf{stationary1}--\rf{stationary2}, we obtain \rf{stationary}.
\endpf

\ms

Note that \rf{FBSDEi} for $i=1,2$ are two coupled FBSDEs, and these two FBSDEs are coupled through the relation \rf{stationary}. In fact, from \rf{stationary}, we see that
$$\begin{pmatrix}R^1_{11}&R^1_{12}\\ R^2_{21}&R^2_{22}\end{pmatrix}\begin{pmatrix}u^*_1\\ u^*_2\end{pmatrix}=-\begin{pmatrix}B_1^\top Y_1^*+D_1^\top Z^*_1+S^1_1X^*+\rho^1_1\\
B_2^\top Y_2^*+D_2^\top Z^*_2+S^2_2X^*+\rho^2_2\end{pmatrix}.$$
Thus, say, in the case that the coefficient matrix of $u^*$ is invertible, one has
$$\begin{pmatrix}u^*_1\\ u^*_2\end{pmatrix}=-\begin{pmatrix}R^1_{11}&R^1_{12}\\ R^2_{21}&R^2_{22}\end{pmatrix}^{-1}\begin{pmatrix}B_1^\top Y_1^*+D_1^\top Z^*_1+S^1_1X^*+\rho^1_1\\
B_2^\top Y_2^*+D_2^\top Z^*_2+S^2_2X^*+\rho^2_2\end{pmatrix}.$$
Plugging the above into \rf{FBSDEi}, we see the coupling between the two coupled FBSDEs (with $i=1,2$).

\ms

To conclude this section, let us write FBSDE \rf{FBSDEi} and stationarity condition \rf{stationary} more compactly.
For this, we introduce the following:
$$\ba{ll}
\no\ms\ds\BA(\cd)=\begin{pmatrix}A(\cd)&0 \\ 0&A(\cd)\end{pmatrix},
       \q\BB(\cd)=\begin{pmatrix}B(\cd)&0 \\ 0&B(\cd)\end{pmatrix}
\equiv\begin{pmatrix}B_1(\cd)&B_2(\cd)&0&0\\ 0&0&B_1(\cd)&B_2(\cd)\end{pmatrix},\\
\no\ms\ds\BC(\cd)=\begin{pmatrix}C(\cd)&0 \\ 0&C(\cd)\end{pmatrix},
       \q\BD(\cd)=\begin{pmatrix}D(\cd)&0 \\ 0&D(\cd)\end{pmatrix}
\equiv\begin{pmatrix}D_1(\cd)&D_2(\cd)&0&0\\ 0&0&D_1(\cd)&D_2(\cd)\end{pmatrix},\\
\no\ms\ds\BQ(\cd)=\begin{pmatrix}Q^1(\cd)&0 \\ 0&Q^2(\cd)\end{pmatrix},
       \q\BS(\cd)=\begin{pmatrix}S^1(\cd)&0 \\ 0&S^2(\cd)\end{pmatrix},
       \q\BR(\cd)=\begin{pmatrix}R^1(\cd)&0 \\ 0&R^2(\cd)\end{pmatrix},\\
\ds\Bq(\cd)=\begin{pmatrix}   q^1(\cd) \\ q^2(\cd)   \end{pmatrix},\q
 \Brho(\cd)=\begin{pmatrix}\rho^1(\cd) \\ \rho^2(\cd)\end{pmatrix},\q
        \BG=\begin{pmatrix}G^1&0 \\ 0&G^2\end{pmatrix},\q
        \Bg=\begin{pmatrix}  g^1 \\ g^2  \end{pmatrix}.\ea$$
Then
$$\left\{\2n\ba{ll}
\ds\BA(\cd)\in L^1(0,T;\dbR^{2n\times2n}),\q\BB(\cd)\in L^2(0,T;\dbR^{2n\times2m}),\\
\ns\ds\BC(\cd)\in L^2(0,T;\dbR^{2n\times2n}),\q\BD(\cd)\in L^\i(0,T;\dbR^{2n\times2m}),\\
\ns\ds\BQ(\cd)\in L^1(0,T;\dbS^{2n}),\q\BS(\cd)\in L^2(0,T;\dbR^{2m\times2n}),
\q\BR(\cd)\in L^\i(0,T;\dbS^{2m}),\\
\ns\ds\Bq(\cd)\in L^2_\dbF(\O;L^1(0,T;\dbR^{2n})),\q\Brho(\cd)\in L_\dbF^2(0,T;\dbR^{2m}),
\q\BG\in\dbS^{2n},\q\Bg\in L^2_{\cF_T}(\O;\dbR^{2n}).\ea\right.$$
Further, let
$$\BJ=\begin{pmatrix}I_{m_1}&0 \\ 0&0 \\ 0&0 \\ 0&I_{m_2}\end{pmatrix}
\equiv\begin{pmatrix}I_{m_1}           & 0_{m_1\times m_2} \\
                     0_{m_2\times m_1} & 0_{m_2\times m_2} \\
                     0_{m_1\times m_1} & 0_{m_1\times m_2} \\
                     0_{m_2\times m_1} & I_{m_2}           \end{pmatrix}\in\dbR^{2m\times m},
\qq\BI_k=\begin{pmatrix}I_k \\ I_k\end{pmatrix}\in\dbR^{2k\times k}.$$
Clearly, one has
$$\ba{ll}
\no\ms\ds\BB(\cd)\BJ
\equiv\begin{pmatrix}B_1(\cd)&B_2(\cd)&0&0 \\ 0&0&B_1(\cd)&B_2(\cd)\end{pmatrix}
      \begin{pmatrix}I_{m_1}&0 \\ 0&0 \\ 0&0 \\ 0&I_{m_2}\end{pmatrix}
=\begin{pmatrix}B_1(\cd)&0 \\ 0&B_2(\cd)\end{pmatrix},\\
\no\ms\ds\BD(\cd)\BJ
\equiv\begin{pmatrix}D_1(\cd)&D_2(\cd)&0&0 \\ 0&0&D_1(\cd)&D_2(\cd)\end{pmatrix}
      \begin{pmatrix}I_{m_1}&0 \\ 0&0 \\ 0&0 \\ 0&I_{m_2}\end{pmatrix}
=\begin{pmatrix}D_1(\cd)&0 \\ 0&D_2(\cd)\end{pmatrix},\\
\no\ms\ds\BJ^\top\BS(\cd)
\equiv\begin{pmatrix}I_{m_1}&0&0&0 \\ 0&0&0&I_{m_2}\end{pmatrix}
      \begin{pmatrix}S^1_1(\cd) & 0 \\ S^1_2(\cd) & 0\\
                     0 & S^2_1(\cd) \\ 0 & S^2_2(\cd)\end{pmatrix}
=\begin{pmatrix}S^1_1(\cd)&0 \\ 0&S^2_2(\cd)\end{pmatrix},\\
\no\ms\ds\BJ^\top\BR(\cd)
\equiv\begin{pmatrix}I_{m_1}&0&0&0 \\ 0&0&0&I_{m_2}\end{pmatrix}
      \begin{pmatrix}R^1_1(\cd)&0 \\ R^1_2(\cd)&0\\
                     0&R^2_1(\cd) \\ 0&R^2_2(\cd)\end{pmatrix}
=\begin{pmatrix}R^1_1(\cd)&0 \\ 0&R^2_2(\cd)\end{pmatrix},\\
\ds\BJ^\top\Brho(\cd)
\equiv\begin{pmatrix}I_{m_1}&0&0&0 \\ 0&0&0&I_{m_2}\end{pmatrix}
      \begin{pmatrix}\rho^1_1(\cd) \\ \rho^1_2(\cd)\\
                     \rho^2_1(\cd) \\ \rho^2_2(\cd)\end{pmatrix}
=\begin{pmatrix}\rho^1_1(\cd) \\ \rho^2_2(\cd)\end{pmatrix}.\ea$$
With the above notation, FBSDE \rf{FBSDEi} can be written as (suppressing $s$ and dropping $*$)
\bel{FBSDE**}\left\{\2n\ba{ll}
\ds dX=\big(AX+Bu+b\big)ds+\big(CX+Du+\si\big)dW,\\
\ns\ds d\BY=-\big(\BA^\top\BY+\BC^\top\BZ+\BQ\BI_nX+\BS^\top\BI_mu+\Bq\big)ds+\BZ dW,\\
\ns\ds X(t)=x,\qq \BY(T)=\BG\BI_nX(T)+\Bg,\ea\right.\ee
where
$$\BY(\cd)=\begin{pmatrix}Y_1(\cd)\\ Y_2(\cd)\end{pmatrix},\qq\BZ(\cd)=\begin{pmatrix}Z_1(\cd)\\ Z_2(\cd)\end{pmatrix}.$$
and the stationarity condition \rf{stationary} can be written as
\bel{stationary**}\BJ^\top\big(\BB^\top\BY+\BD^\top\BZ+\BS\BI_nX+\BR\BI_m u+\Brho\big)=0,\qq\ae~s\in[t,T],~\as\ee
Keep in mind that \rf{FBSDE**} is a coupled FBSDE with the coupling given through \rf{stationary**}.

\section{Closed-Loop Nash Equilibria and Riccati Equations}

We now look at closed-loop Nash equilibria for Problem (SDG). Again, for simplicity of notation,
we will suppress the time variable $s$ as long as no confusion arises. First, we present the
following result which is a consequence of Theorem \ref{Theorem 4.1}.

\bp{prop 5.1}\sl Let {\rm(G1)--(G2)} hold. Suppose that
$(\Th^*_1(\cd),v_1^*(\cd);\Th^*_2(\cd),v^*_2(\cd))
\in\cQ_1[t,T]\times\cU_1[t,T]\times\cQ_2[t,T]\times\cU_2[t,T]$
is a closed-loop Nash equilibrium of Problem {\rm(SDG)} on $[t,T]$.
Denote $\Th^*(\cd)\1n\equiv\1n(\Th^*_1(\cd)\1n^\top\1n,\Th^*_2(\cd)\1n^\top)\1n^\top$
and let $\dbX(\cd)$ be the solution to the $\dbR^{n\times n}$-valued SDE
\bel{16Apr13-dbX}\left\{\2n\ba{ll}
\ds d\dbX=(A+B\Th^*)\dbX ds+(C+D\Th^*)\dbX dW,\q~s\in[t,T],\\
\ns\ds \dbX(t)=I. \ea\right.\ee
Then for $i=1,2$, the adapted solution $(\dbY_i(\cd),\dbZ_i(\cd))$ to the
$\dbR^{n\times n}$-valued BSDE
\bel{16Apr13-dbYi}\left\{\2n\ba{ll}
\ds d\dbY_i=-\Big\{(A+B\Th^*)^\top\dbY_i+(C+D\Th^*)^\top\dbZ_i\\
\ns\ds\qq\qq~
+\big[Q^i+(\Th^*)^\top S^i+(S^i)^\top\Th^*+(\Th^*)^\top R^i\Th^*\big]\dbX\Big\}ds
+\dbZ_idW,\q~s\in[t,T],\\
\ns\ds \dbY_i(T)=G^i\dbX(T),\ea\rt.\ee
satisfies
\bel{16Apr13-yushu}B_i^\top\dbY_i+D_i^\top\dbZ_i+(S_i^i+R_i^i\Th^*)\dbX=0,\q~\ae~s\in[t,T],~\as\ee
\ep

\it Proof. \rm Let us consider state equation \rf{16Apr4-20:30} with the cost functionals
defined by \rf{16Apr4-20:40}. Denoting $v(\cd)=(v_1(\cd)^\top,v_2(\cd)^\top)^\top$,
by an argument similar to \rf{16Apr7-22:00}, we have:
$$\ba{ll}
\ds\wt J^i(t,x;v(\cd))\equiv J^i(t,x;\Th^*(\cd)X(\cd)+v(\cd))\\
\ns\ds=\dbE\Big\{\lan G^iX(T),X(T)\ran+2\lan g^i,X(T)\ran\\
\ns\ds\q+\1n\int_t^T\1n\[\lan{\sc
\begin{pmatrix}\sc Q^i+(\Th^*)^\top\1n S^i+(S^i)^\top\1n\Th^*+(\Th^*)^\top\1n R^i\Th^*
             &\1n\sc (S^i)^\top\1n+(\Th^*)^\top\1n R^i \\ \sc S^i+R^i\Th^* &\1n\sc R^i \end{pmatrix}
\begin{pmatrix}\sc X \\ \sc v\end{pmatrix},
\begin{pmatrix}\sc X \\ \sc v\end{pmatrix}}\ran\1n
+\1n2\lan{\sc\begin{pmatrix}\sc q^i+(\Th^*)^\top\1n\rho^i \\ \sc\rho^i\end{pmatrix},
       \begin{pmatrix}\sc X \\ \sc v \end{pmatrix}}\ran\]ds\Big\}.\ea$$
We know by (ii) of Proposition \ref{bp-16Apr4-17:30} that $v^*(\cd)\equiv(v_1^*(\cd)^\top,v^*_2(\cd)^\top)^\top$
is an open-loop Nash equilibrium for the problem with the state equation \rf{16Apr4-20:30} and with the
cost functionals $\wt J^i(t,x;v(\cd))$ for any initial pair $(t,x)$. Thus, according to Theorem \ref{Theorem 4.1},
we have for $i=1,2$,
\bel{16Apr13-14:00}B_i^\top Y_i^*+D_i^\top Z_i^*+(S_i^i+R_i^i\Th^*)X^*+R^i_iv^*+\rho_i^i=0,
\q~\ae~s\in[t,T],~\as\ee
with $X^*(\cd)$ being the solution to the closed-loop system:
\bel{16Apr13-14:25}\left\{\2n\ba{ll}
\ds dX^*=\big[(A+B\Th^*)X^*+Bv^*+b\big]ds+\big[(C+D\Th^*)X^*+Dv^*+\si\big]dW,\q~s\in[t,T],\\
\ns\ds X^*(t)=x, \ea\right.\ee
and $(Y_i^*(\cd),Z_i^*(\cd))$ being the adapted solution to the following BSDE:
\bel{16Apr13-15:00}\left\{\2n\ba{ll}
\ds dY_i^*=-\Big\{(A\1n+\1nB\Th^*)^\top\1nY_i^*+(C\1n+\1nD\Th^*)^\top\1n Z_i^*
+\big[Q^i\1n+\1n(\Th^*)^\top\1n S^i\1n+\1n(S^i)^\top\1n\Th^*\1n+\1n(\Th^*)^\top\1n R^i\Th^*\big]X^*\\
\ns\ds\qq\qq~+(S^i+R^i\Th^*)^\top v^*+q^i+(\Th^*)^\top\rho^i\Big\}ds+Z_i^*dW,\q~s\in[t,T],\\
\ns\ds Y_i^*(T)=G^iX^*(T)+g^i.\ea\rt.\ee
Since $(\Th^*(\cd),v^*(\cd))$ is independent of $x$ and \rf{16Apr13-14:00}--\rf{16Apr13-15:00} hold
for all $x\in\dbR^n$, by subtracting solutions corresponding to $x$ and $0$, the latter from the
former, we see that for any $x\in\dbR^n$, the adapted solution $(X(\cd),Y_i(\cd),Z_i(\cd))$
$(i=1,2)$ to the following FBSDE:
$$\left\{\2n\ba{ll}
\ds dX=(A+B\Th^*)Xds+(C+D\Th^*)XdW,\q~s\in[t,T],\\
\ns\ds dY_i=-\Big\{(A+B\Th^*)^\top Y_i+(C+D\Th^*)^\top Z_i\\
\ns\ds\qq\qq~
+\big[Q^i+(\Th^*)^\top S^i+(S^i)^\top\Th^*+(\Th^*)^\top R^i\Th^*\big]X\Big\}ds+Z_idW,\q~s\in[t,T],\\
\ns\ds X(t)=x,\qq Y_i(T)=G^iX(T),\ea\rt.$$
satisfies
$$B_i^\top Y_i+D_i^\top Z_i+(S_i^i+R_i^i\Th^*)X=0,\q~\ae~s\in[t,T],~\as$$
The desired result then follows easily.
\endpf

\ms

Now we are ready to present the main result of this section,
which characterizes the closed-loop Nash equilibrium of Problem (SDG).

\bt{Theorem 5.2}\sl Let {\rm(G1)--(G2)} hold. Then $(\Th^*(\cd),v^*(\cd))\in\cQ[t,T]\times\cU[t,T]$
is a closed-loop Nash equilibrium of Problem {\rm(SDG)} on $[t,T]$ if and only if the following hold:

\ms

{\rm(i)} For $i=1,2$, the solution $P_i(\cd)\in C([t,T];\dbS^n)$ to the Lyapunov type equation
\bel{Riccati-i}\left\{\2n\ba{ll}
\ds\dot P_i+P_iA+A^\top P_i+C^\top P_iC+Q^i+(\Th^*)^\top(R^i+D^\top P_iD)\Th^*\\
\ns\ds~~~+\big[P_iB+C^\top P_iD+(S^i)^\top\big]\Th^*+(\Th^*)^\top\big[B^\top P_i+D^\top P_iC+S^i\big]=0,
\q~\ae~s\in[t,T],\\
\ns\ds P_i(T)=G^i,\ea\right.\ee
satisfies the following two conditions:
\begin{eqnarray}
&\label{R+DPD>0}R_{ii}^i+D_i^\top P_iD_i\ges0,\q~\ae~s\in[t,T],&\\
&\label{stationary-closed-i}B_i^\top P_i+D_i^\top P_iC+S_i^i+(R_i^i+D_i^\top P_iD)\Th^*=0,\q~\ae~s\in[t,T].&
\end{eqnarray}

\ms

{\rm(ii)} For $i=1,2$, the adapted solution $(\eta_i(\cd),\z_i(\cd))$ to the BSDE
\bel{BSDE-closed-i}\left\{\2n\ba{ll}
\ds d\eta_i=-\Big\{A^\top\eta_i+C^\top\z_i
+(\Th^*)^\top\big[B^\top\eta_i+D^\top\z_i+D^\top\1n P_i\si+\rho^i+(R^i+D^\top\1n P_iD)v^*\big]\\
\ns\ds\qq\q~~+\big[P_iB+C^\top P_iD+(S^i)^\top\big]v^*+C^\top P_i\si+P_ib+q^i\Big\}ds+\z_i dW,\q~s\in[t,T],\\
\ns\ds \eta_i(T)=g^i, \ea\right.\ee
satisfies
\bel{eta-i}B_i^\top\eta_i+D_i^\top\z_i+D_i^\top P_i\si+\rho_i^i+(R_i^i+D_i^\top P_iD)v^*=0,
\q~\ae~s\in[t,T],~\as\ee
\et

\it Proof. \rm We first prove the necessity.
Suppose that $(\Th^*(\cd),v^*(\cd))$ is a closed-loop Nash equilibrium of Problem (SDG) on $[t,T]$,
where $\Th^*(\cd)\equiv(\Th^*_1(\cd)^\top,\Th^*_2(\cd)^\top)^\top$ and
$v^*(\cd)\equiv(v_1^*(\cd)^\top,v_2^*(\cd)^\top)^\top$. Let $\dbX(\cd)$ and $\dbY_i(\cd)$ $(i=1,2)$
be the solutions of \rf{16Apr13-dbX} and \rf{16Apr13-dbYi}, respectively. Consider the following
linear ordinary differential equation (ODE, for short) which is equivalent to \rf{Riccati-i}:
\bel{Riccati-i*}\left\{\2n\ba{ll}
\ds\dot P_i+P_i(A+B\Th^*)+(A+B\Th^*)^\top P_i+(C+D\Th^*)^\top P_i(C+D\Th^*)\\
\ns\ds~~~+Q^i+(\Th^*)^\top S^i+(S^i)^\top\Th^*+(\Th^*)^\top R^i\Th^*=0,\qq s\in[t,T],\\
\ns\ds P_i(T)=G^i.\ea\right.\ee
Such an equation admits a unique solution $P_i(\cd)\in C([t,T];\dbS^n)$. By It\^o's formula, we have
$$\ba{lll}
\ds d(P_i\dbX)\3n&=&\3n\dot P_i\dbX ds+P_i(A+B\Th^*)\dbX ds+P_i(C+D\Th^*)\dbX dW\\
\ns\3n&=&\3n-\,\big\{(A+B\Th^*)^\top P_i\dbX +(C+D\Th^*)^\top P_i(C+D\Th^*)\dbX\\
\ns\3n&~&\3n\q~+\big[Q^i+(\Th^*)^\top S^i+(S^i)^\top\Th^*+(\Th^*)^\top R^i\Th^*\big]\dbX\big\}ds
+P_i(C+D\Th^*)\dbX dW.  \ea$$
Comparing the above with \rf{16Apr13-dbYi}, by the uniqueness of adapted solutions to BSDEs, one has
$$\dbY_i=P_i\dbX,\q \dbZ_i=P_i(C+D\Th^*)\dbX;\qq i=1,2.$$
From \rf{16Apr13-dbX}, we see that the process $\dbX(\cd)$ is invertible almost surely.
Then, the above together with \rf{16Apr13-yushu} leads to \rf{stationary-closed-i}.
Now let $X^*(\cd)$ be the solution to \rf{16Apr13-14:25}, and for $i=1,2$,
let $(Y_i^*(\cd),Z_i^*(\cd))$ be the adapted solution to \rf{16Apr13-15:00}. Define
\bel{etaz}\left\{\2n\ba{ll}
\ds\eta_i=Y_i^*-P_iX^*,\\
\ns\z_i=Z_i^*-P_i(C+D\Th^*)X^*-P_i(Dv^*+\si).\ea\right.\ee
Then $\eta_i(T)=g^i$, and
$$\ba{lll}
\ds d\eta_i \3n&=&\3n dY_i^*-\dot P_iX^*ds-P_i dX^*\\
\ns\3n&=&\3n -\,\Big\{(A+B\Th^*)^\top Y_i^*+(C+D\Th^*)^\top Z_i^*
+(S^i+R^i\Th^*)^\top v^*+P_i(Bv^*+b)+q^i+(\Th^*)^\top\rho^i\\
\ns\3n&~&\3n\q~
+\big[\dot P_i+P_i(A+B\Th^*)+Q^i+(\Th^*)^\top S^i+(S^i)^\top\Th^*+(\Th^*)^\top R^i\Th^*\big]X^*\Big\}ds\\
\ns\3n&~&\3n +\,\Big\{Z_i^*-P_i\big[(C+D\Th^*)X^*+Dv^*+\si\big]\Big\}dW\\
\ns\3n&=&\3n -\,\Big\{(A+B\Th^*)^\top Y_i^*+(C+D\Th^*)^\top Z_i^*
+(S^i+R^i\Th^*)^\top v^*+P_i(Bv^*+b)+q^i+(\Th^*)^\top\rho^i\\
\ns\3n&~&\3n\q~
-(A+B\Th^*)^\top P_iX^*-(C+D\Th^*)^\top P_i(C+D\Th^*)X^*\Big\}ds+\z_idW\\
\ns\3n&=&\3n -\,\Big\{(A+B\Th^*)^\top\eta_i+(C+D\Th^*)^\top\z_i+(C+D\Th^*)^\top P_i(Dv^*+\si)\\
\ns\3n&~&\3n\q~+(S^i+R^i\Th^*)^\top v^*+P_i(Bv^*+b)+q^i+(\Th^*)^\top\rho^i\Big\}ds+\z_idW\\
\ns\3n&=&\3n -\,\Big\{A^\top\eta_i+C^\top\z_i
+(\Th^*)^\top\big[B^\top\eta_i+D^\top\z_i+D^\top P_i\si+\rho^i+(R^i+D^\top P_iD)v^*\big]\\
\ns\3n&~&\3n\q~+\big[P_iB+C^\top P_iD+(S^i)^\top\big]v^*+C^\top P_i\si+P_ib+q^i\Big\}ds+\z_idW.\ea$$
Thus, $(\eta_i,\z_i)$ is the adapted solution to BSDE \rf{BSDE-closed-i}.
Next, from the proof of Proposition \ref{prop 5.1} we know that \rf{16Apr13-14:00} holds.
Thus (noting \rf{stationary-closed-i} and \rf{etaz}),
$$\ba{lll}
\ds0 \4n&=&\4n\ds B_i^\top Y_i^*+D_i^\top Z_i^*+(S_i^i+R_i^i\Th^*)X^*+R^i_iv^*+\rho_i^i\\
\ns\4n&=&\4n\ds
B_i^\top\eta_i\1n+\1n D_i^\top\z_i\1n+\1n D_i^\top P_i\si\1n+\1n\rho_i^i\1n+\1n(R_i^i\1n+\1n D_i^\top P_iD)v^*\1n
+\1n\big[B_i^\top P_i\1n+\1n D_i^\top P_iC\1n+\1n S_i^i\1n+\1n(R_i^i\1n+\1n D_i^\top P_iD)\Th^*\big]X^*\\
\ns\4n&=&\4n\ds B_i^\top\eta_i+D_i^\top\z_i+D_i^\top P_i\si+\rho_i^i+(R_i^i+D_i^\top P_iD)v^*,\ea$$
which is \rf{eta-i}. The proof of \rf{R+DPD>0} will be included in the proof of sufficiency.

\ms

To prove the sufficiency, we take any $v(\cd)=(v_1(\cd)^\top,v_2(\cd)^\top)^\top\in\cU_1[t,T]\times\cU_2[t,T]$.
Denote $w=(v_1^\top$, $(v_2^*)^\top)^\top$, and let
$$X(\cd)=X(\cd\,;t,x,\Th_1^*(\cd),v_1(\cd),\Th_2^*(\cd),v_2^*(\cd))$$
be the state process corresponding to $(t,x)$ and $(\Th_1^*(\cd),v_1(\cd),\Th_2^*(\cd),v_2^*(\cd))$.
By It\^{o}'s formula, we have
$$\ba{lll}
\ds\dbE\[\blan G^1X(T),X(T)\bran+2\blan g^1,X(T)\bran\]-\dbE\[\lan P_1(t)x,x\ran+2\lan \eta_1(t),x\ran\]\\
\ns\ds=\dbE\int_t^T\Big\{\blan\dot P_1X,X\bran+2\lan P_1X,(A+B\Th^*)X+Bw+b\ran\\
\ns\ds\qq\qq~+\blan P_1\big[(C+D\Th^*)X+Dw+\si\big],(C+D\Th^*)X+Dw+\si\bran\\
\ns\ds\qq\qq~-2\blan A^\top\eta_1+C^\top\z_1
+(\Th^*)^\top\big[B^\top\eta_1+D^\top\z_1+D^\top\1n P_1\si+\rho^1+(R^1+D^\top\1n P_1D)v^*\big],X\bran\\
\ns\ds\qq\qq~-2\blan\big[P_1B+C^\top P_1D+(S^1)^\top\big]v^*+C^\top P_1\si+P_1b+q^1,X\bran\\
\ns\ds\qq\qq~+2\lan\eta_1,(A+B\Th^*)X+Bw+b\ran+2\lan\z_1,(C+D\Th^*)X+Dw+\si\ran\Big\}ds\ea$$
$$\ba{ll}
\ns\ds=\dbE\int_t^T\Big\{\blan\big[\dot P_1+P_1(A+B\Th^*)+(A+B\Th^*)^\top P_1
+(C+D\Th^*)^\top P_1(C+D\Th^*)\big]X,X\bran\\
\ns\ds\qq\qq~+2\lan P_1X,Bw+b\ran+2\lan P_1(C+D\Th^*)X,Dw+\si\ran+\lan P_1(Dw+\si),Dw+\si\ran\\
\ns\ds\qq\qq~-2\blan(\Th^*)^\top\big[D^\top\1n P_1\si+\rho^1+(R^1+D^\top\1n P_1D)v^*\big],X\bran\\
\ns\ds\qq\qq~-2\blan\big[P_1B+C^\top P_1D+(S^1)^\top\big]v^*+C^\top P_1\si+P_1b+q^1,X\bran\\
\ns\ds\qq\qq~+2\lan\eta_1,Bw+b\ran+2\lan\z_1,Dw+\si\ran\Big\}ds\\
\ns\ds=\dbE\int_t^T\Big\{\blan\big[\dot P_1+P_1(A+B\Th^*)+(A+B\Th^*)^\top P_1
+(C+D\Th^*)^\top P_1(C+D\Th^*)\big]X,X\bran\\
\ns\ds\qq\qq~+2\blan(P_1B+C^\top P_1D)w-\big[P_1B+C^\top P_1D+(S^1)^\top\big]v^*-q^1,X\bran\\
\ns\ds\qq\qq~+2\blan D^\top P_1Dw-(R^1+D^\top P_1D)v^*-\rho^1,\Th^*X\bran+\blan D^\top P_1Dw,w\bran\\
\ns\ds\qq\qq~+2\blan B^\top\eta_1+D^\top\z_1+D^\top P_1\si,w\bran
+\lan P_1\si,\si\ran+2\lan\eta_1,b\ran+2\lan\z_1,\si\ran\Big\}ds.\ea$$
On the other hand, we have
$$\ba{ll}
\ds J^1(t,x;\Th^*X(\cd)+w(\cd))-\dbE\[\blan G^1X(T),X(T)\bran+2\blan g^1,X(T)\bran\]\\
\ns\ds=\dbE\int_t^T\[\blan\sc{
\begin{pmatrix}\sc Q^1+(\Th^*)^\top S^1+(S^1)^\top\Th^*+(\Th^*)^\top R^1\Th^*
              &\sc (S^1)^\top+(\Th^*)^\top R^1 \\ \sc S^1+R^1\Th^* & \sc R^1 \end{pmatrix}
\begin{pmatrix}\sc X \\ \sc w\end{pmatrix},
\begin{pmatrix}\sc X \\ \sc w\end{pmatrix}}\bran+2\blan{\sc
\begin{pmatrix}\sc q^1+(\Th^*)^\top\rho^1 \\ \sc\rho^1\end{pmatrix},
\begin{pmatrix}\sc X \\ \sc w \end{pmatrix}}\bran\]ds\\
\ns\ds=\dbE\int_t^T\Big\{\blan\big[Q^1+(\Th^*)^\top S^1+(S^1)^\top\Th^*
+(\Th^*)^\top R^1\Th^*\big]X,X\bran+2\blan(S^1)^\top w+q^1,X\bran\\
\ns\ds\qq\qq\q~+2\blan R^1w+\rho^1,\Th^*X\bran+\blan R^1w,w\bran+2\blan\rho^1,w\bran\Big\}ds.\ea$$
Combining the above two equations, together with equation \rf{Riccati-i*} (which is equivalent to \rf{Riccati-i})
and conditions \rf{stationary-closed-i} and \rf{eta-i}, one obtains
$$\ba{ll}
\ds J^1(t,x;\Th^*X(\cd)+w(\cd))-\dbE\[\lan P_1(t)x,x\ran+2\lan \eta_1(t),x\ran\]\\
\ns\ds=\dbE\int_t^T\Big\{2\blan\big[P_1B+C^\top P_1D+(S^1)^\top\big](w-v^*),X\bran\\
\ns\ds\qq\qq~+2\blan(R^1+D^\top P_1D)(w-v^*),\Th^*X\bran+\blan(R^1+D^\top P_1D)w,w\bran\\
\ns\ds\qq\qq~+2\blan B^\top\eta_1+D^\top\z_1+D^\top P_1\si+\rho^1,w\bran
+\lan P_1\si,\si\ran+2\lan\eta_1,b\ran+2\lan\z_1,\si\ran\Big\}ds\\
\ns\ds=\dbE\int_t^T\Big\{2\blan\big[P_1B_1+C^\top P_1D_1+(S^1_1)^\top\big](v_1-v_1^*),X\bran
+2\blan(R^1_1+D_1^\top P_1D)^\top(v_1-v_1^*),\Th^*X\bran \\
\ns\ds\qq\qq~+\blan(R^1_{11}+D_1^\top P_1D_1)v_1,v_1\bran
+2\blan(R^1_{12}+D_1^\top P_1D_2)v_2^*,v_1\bran+\blan(R^1_{22}+D_2^\top P_1D_2)v_2^*,v_2^*\bran\\
\ns\ds\qq\qq~+2\blan B_1^\top\eta_1+D_1^\top\z_1+D_1^\top P_1\si+\rho_1^1,v_1\bran
+2\blan B_2^\top\eta_1+D_2^\top\z_1+D_2^\top P_1\si+\rho_2^1,v_2^*\bran\\
\ns\ds\qq\qq~+\lan P_1\si,\si\ran+2\lan\eta_1,b\ran+2\lan\z_1,\si\ran\Big\}ds\\
\ns\ds=\dbE\int_t^T\Big\{\blan(R^1_{11}+D_1^\top P_1D_1)v_1,v_1\bran
-2\blan(R^1_{11}+D_1^\top P_1D_1)v_1^*,v_1\bran+\blan(R^1_{22}+D_2^\top P_1D_2)v_2^*,v_2^*\bran\\
\ns\ds\qq\qq~+2\blan B_2^\top\eta_1+D_2^\top\z_1+D_2^\top P_1\si+\rho_2^1,v_2^*\bran
+\lan P_1\si,\si\ran+2\lan\eta_1,b\ran+2\lan\z_1,\si\ran\Big\}ds\\
\ns\ds=\dbE\int_t^T\Big\{\blan(R^1_{11}+D_1^\top P_1D_1)(v_1-v_1^*),v_1-v_1^*\bran
-\blan(R^1_{11}+D_1^\top P_1D_1)v_1^*,v_1^*\bran\\
\ns\ds\qq\qq~+\blan(R^1_{22}+D_2^\top P_1D_2)v_2^*,v_2^*\bran
+2\blan B_2^\top\eta_1+D_2^\top\z_1+D_2^\top P_1\si+\rho_2^1,v_2^*\bran\\
\ns\ds\qq\qq~+\lan P_1\si,\si\ran+2\lan\eta_1,b\ran+2\lan\z_1,\si\ran\Big\}ds.\ea$$
Consequently,
$$ J^1(t,x;\Th^*X(\cd)+w(\cd))-J^1(t,x;\Th^*X^*(\cd)+v^*(\cd))
=\dbE\int_t^T\blan(R^1_{11}+D_1^\top P_1D_1)(v_1-v_1^*),v_1-v_1^*\bran ds.$$
It follows that for any $ v_1(\cd)\in\cU_1[t,T]$,
$$ J^1(t,x;\Th_1^*(\cd)X^*(\cd)+v_1^*(\cd),\Th^*_2(\cd)X(\cd)+v_2^*(\cd))
\les J^1(t,x;\Th_1^*(\cd)X(\cd)+v_1(\cd),\Th^*_2(\cd)X(\cd)+v_2^*(\cd)),$$
if and only if
$$R^1_{11}+D_1^\top P_1D_1\ges0,\qq\ae~s\in[t,T].$$
Similarly, for any $ v_2(\cd)\in\cU_2[t,T]$,
$$ J^2(t,x;\Th_1^*(\cd)X^*(\cd)+v_1^*(\cd),\Th^*_2(\cd)X(\cd)+v_2^*(\cd))
\les J^2(t,x;\Th_1^*(\cd)X(\cd)+v_1^*(\cd),\Th^*_2(\cd)X(\cd)+v_2(\cd)),$$
if and only if
$$R^2_{22}+D_2^\top P_2D_2\ges0,\qq\ae~s\in[t,T].$$
This proves the sufficiency, as well as the necessity of \rf{R+DPD>0}. \endpf

\ms

Note that condition \rf{stationary-closed-i} is equivalent to the following:
$$\begin{pmatrix}B_1^\top P_1+D_1^\top P_1C+S_1^1 \\ B_2^\top P_2+D_2^\top P_2C+S_2^2\end{pmatrix}
 +\begin{pmatrix}R_1^1+D_1^\top P_1D \\ R_2^2+D_2^\top P_2D\end{pmatrix}\Th^*=0.$$
Therefore,
\bel{Th*}\Th^*=-\begin{pmatrix}R_1^1+D_1^\top P_1D \\ R_2^2+D_2^\top P_2D\end{pmatrix}^{-1}
\begin{pmatrix}B_1^\top P_1+D_1^\top P_1C+S_1^1 \\ B_2^\top P_2+D_2^\top P_2C+S_2^2\end{pmatrix},\ee
provided the involved inverse (which is an $\dbR^{m\times m}$-valued function) exists.
By plugging such a $\Th^*(\cd)$ into \rf{Riccati-i}, we see that the equations for
$P_1(\cd)$ and $P_2(\cd)$ are coupled, symmetric, and of Riccati type.

\ms

Now, let us try to rewrite the Riccati equation in a more compact form.
Note that (recalling the notation we introduced in the previous section)
$$\ba{lll}
\ds0 \4n&=&\4n\ds
\begin{pmatrix}B_1^\top P_1+D_1^\top P_1C+S_1^1 \\ B_2^\top P_2+D_2^\top P_2C+S_2^2\end{pmatrix}
+\begin{pmatrix}R_1^1+D_1^\top P_1D \\ R_2^2+D_2^\top P_2D\end{pmatrix}\Th^*\\
\ns\4n&=&\4n\ds
\begin{pmatrix}B_1^\top&0 \\ 0&B_2^\top\end{pmatrix}
\begin{pmatrix}P_1&0 \\ 0&P_2\end{pmatrix}
\begin{pmatrix}I_n \\ I_n\end{pmatrix}
+\begin{pmatrix}D_1^\top&0 \\ 0&D_2^\top\end{pmatrix}
 \begin{pmatrix}P_1&0 \\ 0&P_2\end{pmatrix}
 \begin{pmatrix}  C&0 \\ 0&C  \end{pmatrix}
 \begin{pmatrix}  I_n \\ I_n  \end{pmatrix}\\
\ns\4n&~&\4n\ds
+\begin{pmatrix}I_{m_1}&0&0&0 \\ 0&0&0&I_{m_2}\end{pmatrix}
 \begin{pmatrix}S^1&0 \\ 0&S^2\end{pmatrix}
 \begin{pmatrix}  I_n \\ I_n  \end{pmatrix}\\
\ns\4n&~&\4n\ds
+\lt[\begin{pmatrix}I_{m_1}&0&0&0 \\ 0&0&0&I_{m_2}\end{pmatrix}
     \begin{pmatrix}R^1&0 \\ 0&R^2\end{pmatrix}
     \begin{pmatrix}  I_m \\ I_m\end{pmatrix}
+\begin{pmatrix}D_1^\top&0 \\ 0&D_2^\top\end{pmatrix}
 \begin{pmatrix}     P_1&0 \\ 0&P_2     \end{pmatrix}
 \begin{pmatrix}       D&0 \\ 0&D       \end{pmatrix}
 \begin{pmatrix}       I_m \\ I_m       \end{pmatrix}\rt]\Th^*\\
\ns\4n&\equiv&\4n\ds \BJ^\top\big(\BB^\top\BP+\BD^\top\BP\BC+\BS\big)\BI_n
+\big[\BJ^\top\big(\BR+\BD^\top\BP\BD\big)\BI_m\big]\Th^*,\ea$$
with
$$\BP(\cd)\equiv\begin{pmatrix}P_1(\cd)&0\\ 0&P_2(\cd)\end{pmatrix}.$$
Hence, in the case that $\big[\BJ^\top\big(\BR+\BD^\top\BP\BD\big)\BI_m\big]^{-1}
\equiv\begin{pmatrix}R_1^1+D_1^\top P_1D\\ R_2^2+D_2^\top P_2D\end{pmatrix}^{-1}$
exists and is bounded, we have
\bel{Th**}\Th^*=-\big[\BJ^\top\big(\BR+\BD^\top\BP\BD\big)\BI_m\big]^{-1}
\BJ^\top\big(\BB^\top\BP+\BD^\top\BP\BC+\BS\big)\BI_n,\ee
which is the same as \rf{Th*}. On the other hand, \rf{Riccati-i} can be written as
$$\ba{ll}
\ds0=\begin{pmatrix}\dot P_1&0 \\ 0&\dot P_2\end{pmatrix}
    +\begin{pmatrix}P_1&0 \\ 0&P_2\end{pmatrix}
     \begin{pmatrix}A&0 \\ 0&A\end{pmatrix}
    +\begin{pmatrix}A&0 \\ 0&A\end{pmatrix}^\top
     \begin{pmatrix}P_1&0 \\ 0&P_2\end{pmatrix}\ea$$
$$\ba{ll}
\ns\ds\qq+\begin{pmatrix}  C&0 \\ 0&C\end{pmatrix}^\top
          \begin{pmatrix}P_1&0 \\ 0&P_2\end{pmatrix}
          \begin{pmatrix}C&0 \\ 0&C\end{pmatrix}
+\begin{pmatrix}Q_1&0 \\ 0&Q_2\end{pmatrix}\\
\ns\ds\qq+\begin{pmatrix}\Th^*&0 \\ 0&\Th^*\end{pmatrix}^\top
          \begin{pmatrix}R^1+D^\top P_1D&0 \\ 0&R^2+D^\top P_2D\end{pmatrix}
          \begin{pmatrix}\Th^*&0 \\ 0&\Th^*\end{pmatrix}\\
\ns\ds\qq+\begin{pmatrix}P_1B+C^\top P_1D+(S^1)^\top&0 \\ 0&P_2B+C^\top P_2D+(S^2)^\top\end{pmatrix}
          \begin{pmatrix}\Th^*&0 \\ 0&\Th^*\end{pmatrix}\\
\ns\ds\qq+\begin{pmatrix}\Th^*&0 \\ 0&\Th^*\end{pmatrix}^\top
          \begin{pmatrix}B^\top P_1+D^\top P_1C+S^1&0 \\ 0&B^\top P_2+D^\top P_2C+S^2\end{pmatrix}.\ea$$
Consequently, one sees that the following holds:
\bel{Riccati**}\left\{\2n\ba{ll}
\ds\dot\BP+\BP\BA+\BA^\top\BP+\BC^\top\BP\BC+\BQ+\BTh^\top\big(\BR+\BD^\top\BP\BD\big)\BTh\\
\ns\ds\,~~+\big(\BP\BB+\BC^\top\BP\BD+\BS^\top\big)\BTh+\BTh^\top\big(\BB^\top\BP+\BD^\top\BP\BC+\BS\big)=0,
\q~\ae~s\in[t,T],\\
\ns\ds\BP(T)=\BG,\ea\right.\ee
where
$$\BTh(\cd)=\begin{pmatrix}\Th^*(\cd)&0\\ 0&\Th^*(\cd)\end{pmatrix},$$
and $\Th^*$ is given by \rf{Th**}. Clearly, \rf{Riccati**} is symmetric.

\section{Two Examples}

From the previous sections, we see that the existence of an open-loop Nash equilibrium is equivalent to
the solvability of a coupled system of two FBSDEs, together with the convexity condition for the cost
functionals (see \rf{convexity}); and that the existence of a closed-loop Nash equilibrium is equivalent
to the solvability of a coupled system of two symmetric Riccati equations satisfying certain type of
non-negativity condition (see \rf{R+DPD>0}). Then a natural question is: Are open-loop and closed-loop
Nash equilibria really different? In this section, we will present two examples showing that they are
indeed different.

\ms

The following example shows that Problem (SDG) may have only open-loop Nash equilibria.

\bex{Example 6.1}\rm Consider the following Problem (SDG) with one-dimensional state equation
$$\left\{\2n\ba{ll}
\ds dX(s)=\big[u_1(s)+u_2(s)\big]ds+\big[u_1(s)-u_2(s)\big]dW(s),\qq s\in[t,1], \\
\ns\ds X(t)=x,\ea\right.$$
and cost functionals
$$ J^1(t,x;u_1(\cd),u_2(\cd))=J^2(t,x;u_1(\cd),u_2(\cd))=\dbE X(1)^2\equiv J(t,x;u_1(\cd),u_2(\cd)).$$
Let $\b\ges{1\over 1-t}$. We claim that
$$\lt(u_1^\b(s),u_2^\b(s)\rt)=-\lt({\b x\over2}{\bf1}_{[t,t+{1\over\b}]}(s),
{\b x\over2}{\bf1}_{[t,t+{1\over\b}]}(s)\rt), \qq s\in[t,1],$$
is an open-loop Nash equilibrium of the problem for the initial pair $(t,x)$.
Indeed, it is clear that for any $u_1(\cd)\in L^2_\dbF(t,1;\dbR)$,
$$J(t,x;u_1(\cd),u_2^\b(\cd))\ges0.$$
On the other hand, the state process $X^\b(\cd)$ corresponding to $\big(u_1^\b(s),u_2^\b(s)\big)$ and $(t,x)$
satisfies $X^\b(1)=0$. Hence,
$$J(t,x;u_1^\b(\cd),u_2^\b(\cd))=0\les J(t,x;u_1(\cd),u_2^\b(\cd)),\qq \forall\, u_1(\cd)\in L^2_\dbF(t,1;\dbR).$$
Likewise,
$$J(t,x;u_1^\b(\cd),u_2^\b(\cd))=0\les J(t,x;u_1^\b(\cd),u_2(\cd)),\qq \forall\, u_2(\cd)\in L^2_\dbF(t,1;\dbR).$$
This establishes the claim.

\ms

However, this problem does {\it not} admit a closed-loop Nash equilibrium. We now show this by contradiction.
Suppose $(\Th^*_1(\cd),v_1^*(\cd);\Th^*_2(\cd),v^*_2(\cd))$ is a closed-loop Nash equilibrium.
Consider the corresponding ODEs in Theorem \ref{Theorem 5.2}, which now become
\bel{16Apr27-20:00}\left\{\2n\ba{ll}
\ds\dot P_i+P_i(\Th_1^*-\Th^*_2)^2+2P_i(\Th_1^*+\Th^*_2)=0,\\
\ns\ds P_i(1)=1,\ea\right.\qq i=1,2.\ee
The corresponding constraints read
\bel{6.4}P_1,P_2\ges0,\qq P_1+P_1(\Th_1^*-\Th^*_2)=0,\qq
P_2-P_2(\Th_1^*-\Th^*_2)=0.\ee
Since $P_1(\cd)$ and $P_2(\cd)$ satisfy the same ODE \rf{16Apr27-20:00}, we have $P_1(\cd)=P_2(\cd)$.
Then \rf{6.4} implies $P_1(\cd)=0$, which contradicts the terminal condition $P_1(1)=1$.
\ex

The following example shows that Problem (SDG) may have only closed-loop Nash equilibria.

\bex{Example 6.2}\rm Consider the following Problem (SDG) with one-dimensional state equation
$$\left\{\2n\ba{ll}
\ds dX(s)=u_1(s)ds+u_2(s)dW(s),\qq s\in[t,1], \\
\ns\ds X(t)=x,\ea\right.$$
and cost functionals
\begin{eqnarray}
J^1(t,x;u_1(\cd),u_2(\cd))\3n&=&\3n\dbE\Big\{|X(1)|^2+\int_t^1|u_1(s)|^2 ds\Big\},\nonumber\\
J^2(t,x;u_1(\cd),u_2(\cd))\3n&=&\3n\dbE\Big\{-|X(1)|^2+\int_t^1\[-|X(s)|^2+|u_2(s)|^2\] ds\Big\}.\nonumber
\end{eqnarray}
We claim that the problem admits a closed-loop Nash equilibrium of form $(\Th_1(\cd),0;$ $\Th_2(\cd),0)$.
In fact, by Theorem \ref{Theorem 5.2}, we need to solve the following Riccati equations for $P_1(\cd)$
and $P_2(\cd)$:
\begin{eqnarray}
&&\label{16Apr26-17:00}\left\{\2n\ba{ll}
\ds \dot P_1(s)+P_1(s)\Th_2(s)^2+2P_1(s)\Th_1(s)+\Th_1(s)^2=0, \\
\ns\ds P_1(1)=1,\\
\ns\ds P_1(s)+\Th_1(s)=0, \ea\right.\\
\ns&&\label{16Apr26-17:30}\left\{\2n\ba{ll}
\ds \dot P_2(s)+P_2(s)\Th_2(s)^2+2P_2(s)\Th_1(s)+\Th_2(s)^2-1=0, \\
\ns\ds P_2(1)=-1,\\
\ns\ds 1+P_2(s)\ges0,\\
\ns\ds [1+P_2(s)]\Th_2(s)=0.\ea\right.
\end{eqnarray}
By the fourth equation in \rf{16Apr26-17:30}, we may assume $\Th_2(\cd)=0$.
Then \rf{16Apr26-17:00}--\rf{16Apr26-17:30} become
(taking into account $\Th_1(\cd)=-P_1(\cd)$ from the third equation in \rf{16Apr26-17:00})
$$\left\{\2n\ba{ll}
\ds \dot P_1(s)=P_1(s)^2, \\
\ns\ds P_1(1)=1,\ea\right.\qq
\left\{\2n\ba{ll}
\ds \dot P_2(s)=2P_1(s)P_2(s)+1, \\
\ns\ds P_2(1)=-1,\q 1+P_2(s)\ges0.\ea\right. $$
A straightforward calculation leads to
%
%
$$P_1(s)={1\over 2-s},\qq P_2(s)={-(2-s)^3-2\over 3(2-s)^2}.$$
Therefore, $((2-s)^{-1},0;0,0)$ is a closed-loop Nash equilibrium of the problem.

\ms

Next, we claim that the problem does {\it not} have open-loop Nash equilibria.
Indeed, suppose $(u_1^*(\cd),u_2^*(\cd))$ is an open-loop Nash equilibrium for some initial pair $(t,x)$.
Then $u_2^*(\cd)$ is an open-loop optimal control of the following Problem (SLQ) with state equation
\bel{16Apr26-21:00}\left\{\2n\ba{ll}
\ds dX(s)=u_1^*(s)ds+u_2(s)dW(s),\qq s\in[t,1], \\
\ns\ds X(t)=x,\ea\right.\ee
and cost functional
\bel{16Apr26-21:30}\wt J(t,x;u_2(\cd))=\dbE\Big\{-|X(1)|^2+\int_t^1\[-|X(s)|^2+|u_2(s)|^2\] ds\Big\}.\ee
For any $u_2(\cd)\in L^2_\dbF(t,1;\dbR)$, the corresponding solution to \rf{16Apr26-21:00} is given by
\bel{16Apr26-22:00}X(s)=x+\int_t^s u_1^*(r)dr+\int_t^s u_2(r)dW(r).\ee
Let $\e>0$ be undetermined. Substituting \rf{16Apr26-22:00} into \rf{16Apr26-21:30} and using the
inequality $(a+b)^2\ges(1-{1\over\e})a^2+(1-\e)b^2$, we see
$$\ba{lll}
\ds\wt J(t,x;u_2(\cd))\4n&\les&\4n\ds\({1\over\e}-1\)\dbE\(x+\int_t^1u_1^*(s)ds\)^2
+(\e-1)\dbE\(\int_t^1u_2(s)dW(s)\)^2\\
\ns\4n&~&\4n\ds+\,\({1\over\e}-1\)\dbE\int_t^1\(x+\int_t^s u_1^*(r)dr\)^2ds\\
\ns\4n&~&\4n\ds+\,(\e-1)\dbE\int_t^1\(\int_t^s u_2(r)dW(r)\)^2ds+\dbE\int_t^1|u_2(s)|^2 ds\\
\ns\4n&=&\4n\ds\({1\over\e}-1\)\dbE\[\(x+\int_t^1u_1^*(s)ds\)^2
+\int_t^1\(x+\int_t^s u_1^*(r)dr\)^2ds\]\\
\ns\4n&~&\4n\ds+\,\e\dbE\int_t^1|u_2(s)|^2 ds+(\e-1)\dbE\int_t^1\int_t^s|u_2(r)|^2drds.
\ea$$
Now, by taking $u_2(s)=\l$, $\l\in\dbR$, we have
$$\ba{lll}
\ds\wt J(t,x;\l)\4n&\les&\4n\ds\({1\over\e}-1\)\dbE\[\(x+\int_t^1u_1^*(s)ds\)^2
+\int_t^1\(x+\int_t^s u_1^*(r)dr\)^2ds\]\\
\ns\4n&~&\4n\ds+\,{\l^2(1-t)\over2}\big[2\e+(\e-1)(1-t)\big].\ea$$
Choosing $\e>0$ small enough so that $2\e+(\e-1)(1-t)<0$ and then letting $\l\to\i$, we see that
$$\inf_{u_2(\cd)\in L^2_\dbF(t,1;\dbR)}\wt J(t,x;u_2(\cd))=-\i,$$
which contradicts the fact that $u_2^*(\cd)$ is an open-loop optimal control of the associated LQ problem.
\ex

\section{Closed-Loop Representation of Open-Loop Nash Equilibria}

Inspired by the decoupling technique introduced in \cite{Ma-Protter-Yong 1994,Ma-Yong 1999,Yong 1999, Yong 2006},
we now look at the solvability of FBSDE \rf{FBSDEi}--\rf{stationary}. Recall that with the notation introduced
in Section 4, \rf{FBSDEi} and \rf{stationary} are equivalent to \rf{FBSDE**} and \rf{stationary**}, respectively.
To solve FBSDE \rf{FBSDE**}--\rf{stationary**}, let $(\BBeta(\cd),\BBzeta(\cd))$ be the adapted solution to the
following BSDE for some undetermined $\a:[t,T]\times\O\to\dbR^{2n}$:
$$\left\{\2n\ba{ll}
\ds d\BBeta(s)=\a(s)ds+\BBzeta(s)dW(s),\qq s\in[t,T],\\
\ns\ds\BBeta(T)=\Bg,\ea\right.$$
where
$$\BBeta(\cd)=\begin{pmatrix}\eta_1(\cd)\\ \eta_2(\cd)\end{pmatrix},
\qq\BBzeta(\cd)=\begin{pmatrix}\z_1(\cd)\\ \z_2(\cd)\end{pmatrix}.$$
Let $(X(\cd),\BY(\cd),\BZ(\cd))$ be an adapted solution to FBSDE \rf{FBSDE**}. Suppose the following holds:
\bel{BY}\BY(\cd)=\begin{pmatrix}\Pi_1(\cd)X(\cd)+\eta_1(\cd)\\ \Pi_2(\cd)X(\cd)+\eta_2(\cd)\end{pmatrix}\equiv\BPi(\cd) X(\cd)+\BBeta(\cd),\qq\BPi(\cd)\deq\begin{pmatrix}\Pi_1(\cd)\\ \Pi_2(\cd)\end{pmatrix},\ee
for some differentiable maps $\Pi_i:[t,T]\to\dbR^{n\times n}$ with $\Pi_i(T)=G^i$. By It\^o's formula, we have
$$\ba{ll}
\ds-\big(\BA^\top\BY+\BC^\top\BZ+\BQ\BI_nX+\BS^\top\BI_mu+\Bq\big)ds+\BZ dW(s)=d\BY\\
\ns\ds=\big[\dot\BPi X+\BPi(AX+Bu+b)+\a\big]ds+\big[\BPi(CX+Du+\si)+\BBzeta\big]dW(s)\\
\ns\ds=\big[(\dot\BPi+\BPi A)X+\BPi Bu+\BPi b+\a\big]ds+\big[\BPi CX+\BPi Du+\BPi\si+\BBzeta\big]dW(s).\ea$$
Hence, one should have
\bel{BZ}\BZ=\BPi CX+\BPi Du+\BPi\si+\BBzeta.\ee
Then the stationarity condition \rf{stationary**} becomes
$$\ba{lll}
0\4n&=&\4n\ds\BJ^\top\big(\BB^\top\BY+\BD^\top\BZ+\BS\BI_nX+\BR\BI_mu+\Brho\big)\\
\ns\4n&=&\4n\ds\BJ^\top\big[\BB^\top(\BPi X+\BBeta)+\BD^\top(\BPi CX+\BPi Du+\BPi\si
+\BBzeta)+\BS\BI_nX+\BR\BI_mu+\Brho\big]\\
\ns\4n&=&\4n\ds\BJ^\top\big(\BB^\top\BPi+\BD^\top\BPi C+\BS\BI_n\big)X+\BJ^\top\big(\BR\BI_m+\BD^\top\BPi D\big)u+\BJ^\top\big(\BB^\top\BBeta+\BD^\top\BBzeta+\BD^\top\BPi\si+\Brho\big).\ea$$
Note that
$$ \BJ^\top\big(\BR\BI_m+\BD^\top\BPi D\big)
=\begin{pmatrix}   I_{m_1}&0&0&0 \\ 0&0&0&I_{m_2}   \end{pmatrix}
 \begin{pmatrix}R^1+D^\top\Pi_1D \\ R^2+D^\top\Pi_2D\end{pmatrix}
=\begin{pmatrix}R^1_{11}+D_1^\top\Pi_1D_1 & R^1_{12}+D_1^\top\Pi_1D_2\\
                R^2_{21}+D_2^\top\Pi_2D_1 & R^2_{22}+D_2^\top\Pi_2D_2\end{pmatrix}.$$
This is an $\dbR^{m\times m}$-valued function which is not symmetric in general, even $\Pi_1$ and $\Pi_2$ are symmetric.
We now assume that the above is invertible. Then one has
\bel{open-u}\ba{lll}
\ds u\4n&=&\4n\ds-\big[\BJ^\top(\BR\BI_m+\BD^\top\BPi D)\big]^{-1}
\BJ^\top\big(\BB^\top\BPi+\BD^\top\BPi C+\BS\BI_n\big)X\\
\ns\4n&~&\4n\ds-\big[\BJ^\top(\BR\BI_m+\BD^\top\BPi D)\big]^{-1}
\BJ^\top\big(\BB^\top\BBeta+\BD^\top\BBzeta+\BD^\top\BPi\si+\Brho\big),\ea\ee
and
$$\ba{lll}
0\4n&=&\4n\ds\big(\dot\BPi+\BPi A\big)X+\BPi Bu+\BPi b+\a+\BA^\top(\BPi X+\BBeta)\\
\ns\4n&~&\4n\ds+\,\BC^\top(\BPi CX+\BPi Du+\BPi\si+\BBzeta)+\BQ\BI_nX+\BS^\top\BI_mu+\Bq\\
\ns\4n&=&\4n\ds\big(\dot\BPi+\BPi A+\BA^\top\BPi+\BC^\top\BPi C+\BQ\BI_n\big)X
+\big(\BPi B+\BC^\top\BPi D+\BS^\top\BI_m\big)u\\
\ns\4n&~&\4n\ds+\,\a+\BA^\top\BBeta+\BC^\top\BBzeta+\BPi b+\BC^\top\BPi\si+\Bq\\
\ns\4n&=&\4n\ds\big(\dot\BPi+\BPi A+\BA^\top\BPi+\BC^\top\BPi C+\BQ\BI_n\big)X\\
\ns\4n&~&\4n\ds-\,\big(\BPi B+\BC^\top\BPi D+\BS^\top\BI_m\big)\big[\BJ^\top(\BR\BI_m+\BD^\top\BPi D)\big]^{-1}
\BJ^\top\big(\BB^\top\BPi+\BD^\top\BPi C+\BS\BI_n\big)X\\
\ns\4n&~&\4n\ds-\,\big(\BPi B+\BC^\top\BPi D+\BS^\top\BI_m\big)\big[\BJ^\top(\BR\BI_m+\BD^\top\BPi D)\big]^{-1}\BJ^\top\big(\BB^\top\BBeta+\BD^\top\BBzeta+\BD^\top\BPi\si+\Brho\big)\\
\ns\4n&~&\4n\ds+\,\a+\BA^\top\BBeta+\BC^\top\BBzeta+\BPi b+\BC^\top\BPi\si+\Bq \\
\ns\4n&=&\4n\ds\Big\{\dot\BPi+\BPi A+\BA^\top\BPi+\BC^\top\BPi C+\BQ\BI_n\\
\ns\4n&~&\4n\ds-\,\big(\BPi B+\BC^\top\BPi D+\BS^\top\BI_m\big)\big[\BJ^\top(\BR\BI_m+\BD^\top\BPi D)\big]^{-1}
\BJ^\top\big(\BB^\top\BPi+\BD^\top\BPi C+\BS\BI_n\big)\Big\}X\\
\ns\4n&~&\4n\ds-\,\big(\BPi B+\BC^\top\BPi D+\BS^\top\BI_m\big)\big[\BJ^\top(\BR\BI_m+\BD^\top\BPi D)\big]^{-1}
\BJ^\top\big(\BB^\top\BBeta+\BD^\top\BBzeta+\BD^\top\BPi\si+\Brho\big)\\
\ns\4n&~&\4n\ds+\,\a+\BA^\top\BBeta+\BC^\top\BBzeta+\BPi b+\BC^\top\BPi\si+\Bq.\ea$$
Now, let $\BPi(\cd)$ be the solution to the following Riccati equation:
\bel{Riccati-Pi}\left\{\2n\ba{ll}
\ds\dot\BPi+\BPi A+\BA^\top\BPi+\BC^\top\BPi C+\BQ\BI_n\\
\ns\ds-\,\big(\BPi B+\BC^\top\BPi D+\BS^\top\BI_m\big)\big[\BJ^\top(\BR\BI_m\1n+\BD^\top\BPi D)\big]^{-1}
\BJ^\top\big(\BB^\top\BPi+\BD^\top\BPi C+\BS\BI_n\big)=0,\q s\in[t,T],\\
\ns\ds\BPi(T)=\BG\BI_n.\ea\right.\ee
Then the above leads to the BSDE for $(\BBeta(\cd),\BBzeta(\cd))$ of the following form:
\bel{BSDE-open}\left\{\2n\ba{ll}
\ds d\BBeta=-\Big\{\(\BA^\top-\big(\BPi B+\BC^\top\BPi D+\BS^\top\BI_m\big)
\big[\BJ^\top(\BR\BI_m+\BD^\top\BPi D)\big]^{-1}\BJ^\top\BB^\top\)\BBeta\\
\ns\ds\qq\q+\,\(\BC^\top-\big(\BPi B+\BC^\top\BPi D+\BS^\top\BI_m\big)
\big[\BJ^\top(\BR\BI_m+\BD^\top\BPi D)\big]^{-1}\BJ^\top\BD^\top\)\BBzeta\\
\ns\ds\qq\q+\,\(\BC^\top-\big(\BPi B+\BC^\top\BPi D+\BS^\top\BI_m\big)
\big[\BJ^\top(\BR\BI_m+\BD^\top\BPi D)\big]^{-1}\BJ^\top\BD^\top\)\BPi\si\\
\ns\ds\qq\q+\,\BPi b+\Bq-\big(\BPi B+\BC^\top\BPi D+\BS^\top\BI_m\big)
\big[\BJ^\top(\BR\BI_m+\BD^\top\BPi D)\big]^{-1}\BJ^\top\Brho\Big\}ds+\BBzeta dW,\q s\in[t,T],\\
\ns\ds\BBeta(T)=\Bg.\ea\right.\ee
Hence, we have the following result.

\bt{Theorem 4.2} \sl Let {\rm(G1)--(G2)} hold and let $t\in[0,T)$ be given.
Suppose that the convexity condition \rf{convexity} holds for $i=1,2$, and that the Riccati
equation \rf{Riccati-Pi} admits a solution $\BPi(\cd)$. Let $(\BBeta(\cd),\BBzeta(\cd))$ be the
adapted solution to BSDE \rf{BSDE-open} and let $X(\cd)$ be the solution to the following FSDE
with an arbitrary initial state $x$:
\bel{open-closed}\left\{\2n\ba{ll}
\ds dX=\Big\{\(A-B\big[\BJ^\top(\BR\BI_m+\BD^\top\BPi D)\big]^{-1}\BJ^\top
\big(\BB^\top\BPi+\BD^\top\BPi C+\BS\BI_n\big)\) X\\
\ns\ds\qq\qq~-B\big[\BJ^\top(\BR\BI_m+\BD^\top\BPi D)\big]^{-1}\BJ^\top
\big(\BB^\top\BBeta+\BD^\top\BBzeta+\BD^\top\BPi\si+\Brho\big)+b\Big\}ds\\
\ns\ds\qq\q+\,\Big\{\(C-D\big[\BJ^\top(\BR\BI_m+\BD^\top\BPi D)\big]^{-1}\BJ^\top
\big(\BB^\top\BPi+\BD^\top\BPi C+\BS\BI_n\big)\)X\\
\ns\ds\qq\qq~-D\big[\BJ^\top(\BR\BI_m+\BD^\top\BPi D)\big]^{-1}\BJ^\top\big(\BB^\top\BBeta+\BD^\top\BBzeta
+\BD^\top\BPi\si+\Brho\big)+\si\Big\}dW,\q s\in[t,T],\\
\ns\ds X(t)=x.\ea\right.\ee
Then the process $u(\cd)$ defined by \rf{open-u} is an open-loop Nash equilibrium of Problem {\rm(SDG)} for $(t,x)$.

\et

\it Proof. \rm From the above procedure, we see that with $u(\cd)$ defined by \rf{open-u},
the triple $(X(\cd),\BY(\cd)$, $\BZ(\cd))$ defined through FSDE \rf{open-closed}, \rf{BY} and \rf{BZ},
is an adapted solution to FBSDE \rf{FBSDE**}, and that the stationarity condition \rf{stationary**} holds.
Hence, together with the convexity condition \rf{convexity}, making use of Theorem \ref{Theorem 4.1},
we see that $u(\cd)$ is an open-loop Nash equilibrium of Problem (SDG) for $(t,x)$.
\endpf

\ms

Under the assumptions of Theorem \ref{Theorem 4.2}, Problem (SDG) admits an open-loop Nash equilibrium
for every initial state $x$, and the open-loop Nash equilibria take the following form:
\bel{16July9-21:00}u(\cd)=\Th(\cd)X(\cd)+v(\cd),\ee
for some $(\Th(\cd),v(\cd))\in\cQ[t,T]\times\cU[t,T]$ which is independent of $x$.
The above \rf{16July9-21:00} is called a {\it closed-loop representation} of the open-loop
Nash equilibria of Problem (SDG). More precisely, we have the following definition.

\bde{bde-16June29-17:00} \rm We say that open-loop Nash equilibria of Problem {\rm(SDG)}
on $[t,T]$ admit a {\it closed-loop representation}, if there exists a pair $(\Th(\cd),v(\cd))\in
\cQ[t,T]\times\cU[t,T]$ such that for any initial state $x\in\dbR^n$, the process
\bel{16June29-17:40}u(s)\deq \Th(s)X(s)+v(s),\qq s\in[t,T]\ee
is an open-loop Nash equilibrium of Problem {\rm(SDG)} for $(t,x)$,
where $X(\cd)$ is the solution to the following closed-loop system:
\bel{16June29-X}\left\{\2n\ba{ll}
\ds dX(s)=\big\{[A(s)+B(s)\Th(s)]X(s)+B(s)v(s)+b(s)\big\}ds\\
\ns\ds\qq\qq\1n~+\big\{[C(s)+D(s)\Th(s)]X(s)+D(s)v(s)+\si(s)\big\}dW(s),\qq s\in[t,T],\\
\ns\ds X(t)= x. \ea\right.\ee
\ede

Comparing Definitions \ref{bde-16Apr4-17:00} and \ref{bde-16June29-17:00}, it is natural to ask
whether the closed-loop representation of open-loop Nash equilibria is the outcome of some
closed-loop Nash equilibrium. The following example shows that this is not the case in general.

\bex{Example 6.3} \rm Consider the following state equation:
$$\left\{\2n\ba{ll}
\ds dX(s)=\big[u_1(s)+u_2(s)\big]ds+X(s)dW(s),\qq s\in[t,T],\\
\ns\ds X(t)=x,\ea\right.$$
with cost functionals
$$\ba{ll}
\ds J^1(t,x;u_1(\cd),u_2(\cd))=\dbE\[X(T)^2+\int_t^Tu_1(s)^2ds\],\\
\ns\ds J^2(t,x;u_1(\cd),u_2(\cd))=\dbE\[X(T)^2+\int_t^Tu_2(s)^2ds\].\ea$$
For this case, we have
$$\left\{\2n\ba{llll}
\ds A=0,~C=1,     & B_1=B_2=1, & D_1=D_2=0, & b=\si=0,\\
\ns Q^1=Q^2=0, & S^1=S^2=0, & R^1=\begin{pmatrix}1&0 \\ 0&0\end{pmatrix},
                            & R^2=\begin{pmatrix}0&0 \\ 0&1\end{pmatrix},\\
\ns G^1=G^2=1, & q^1=q^2=0, & \rho^1=\rho^2=0, & g^1=g^2=0.\ea\right.$$
Clearly, the convexity condition \rf{convexity} holds for $i=1,2$. In this example,
the Riccati equation \rf{Riccati-Pi} can be written componentwise as follows:
\begin{eqnarray}
&&\label{Pi1-1*}\left\{\2n\ba{ll}
\ds\dot\Pi_1(s)+\Pi_1(s)-\Pi_1(s)\big[\Pi_1(s)+\Pi_2(s)\big]=0,\qq s\in[t,T],\\
\ns\ds\Pi_1(T)=1,\ea\right.\\
\ns&&\label{Pi2-1*}\left\{\2n\ba{ll}
\ds\dot\Pi_2(s)+\Pi_2(s)-\Pi_2(s)\big[\Pi_1(s)+\Pi_2(s)\big]=0,\qq s\in[t,T],\\
\ns\ds\Pi_2(T)=1.\ea\right.
\end{eqnarray}
It is easy to see that
$$\Pi_1(s)=\Pi_2(s)={e^{T-s}\over 2e^{T-s}-1}$$
are solutions to \rf{Pi1-1*} and \rf{Pi2-1*}, respectively. Note that in this case the adapted solution
$(\BBeta(\cd),\BBzeta(\cd))$ to BSDE \rf{BSDE-open} is $(0,0)$. Then by Theorem \ref{Theorem 4.2},
the open-loop Nash equilibria of this Problem (SDG) on $[t,T]$ admit a closed-loop representation given by
\bel{16June26_22:30-1}u_1(s)=u_2(s)=-{e^{T-s}\over 2e^{T-s}-1}X(s),\qq s\in[t,T].\ee

\ss

Next we verify that the problem admits a closed-loop Nash equilibrium of form $(\Th_1(\cd),0;$ $\Th_2(\cd),0)$.
In light of Theorem \ref{Theorem 5.2}, we need to solve the following Riccati equations for $P_1(\cd)$
and $P_2(\cd)$:
\bel{16June26-21:00}\left\{\2n\ba{ll}
\ds\dot P_1(s)+P_1(s)+\Th_1(s)^2+2P_1(s)\big[\Th_1(s)+\Th_2(s)\big]=0,\\
\ns\ds P_1(T)=1,\\
\ns\ds P_1(s)+\Th_1(s)=0, \ea\right.\ee
and
\bel{16June26-21:10}\left\{\2n\ba{ll}
\ds\dot P_2(s)+P_2(s)+\Th_2(s)^2+2P_2(s)\big[\Th_1(s)+\Th_2(s)\big]=0,\\
\ns\ds P_2(T)=1,\\
\ns\ds P_2(s)+\Th_2(s)=0. \ea\right.\ee
Noting the third equations in \rf{16June26-21:00} and \rf{16June26-21:10}, we can further write
\rf{16June26-21:00}-\rf{16June26-21:10} as follows:
\begin{eqnarray}
&&\label{16June26-22:00}\left\{\2n\ba{ll}
\ds\dot P_1(s)=P_1(s)^2+2P_1(s)P_2(s)-P_1(s),\\
\ns\ds P_1(T)=1, \ea\right.\\
\ns&&\label{16June26-22:10}\left\{\2n\ba{ll}
\ds\dot P_2(s)=P_2(s)^2+2P_2(s)P_1(s)-P_2(s),\\
\ns\ds P_2(T)=1. \ea\right.
\end{eqnarray}
Now it is easily seen that
$$P_1(s)=P_2(s)={e^{T-s}\over3e^{T-s}-2}.$$
Hence,
\bel{16June26_22:30-2}\Th_1(s)=\Th_2(s)=-P_1(s)=-{e^{T-s}\over3e^{T-s}-2}.\ee
Comparing \rf{16June26_22:30-1} with \rf{16June26_22:30-2}, we see that the closed-loop representation
of open-loop Nash equilibria is different from the outcome of closed-loop Nash equilibria.
\ex

Now we give a characterization of the closed-loop representation of open-loop Nash equilibria.

\bt{bt-16June29-17:30}\sl Let {\rm(G1)--(G2)} hold and let $(\Th(\cd),v(\cd))\in\cQ[t,T]\times\cU[t,T]$.
Then open-loop Nash equilibria of Problem {\rm(SDG)} on $[t,T]$ admit the closed-loop representation
\rf{16June29-17:40} if and only if the following hold:

\ms

{\rm(i)} The convexity condition \rf{convexity} holds for $i=1,2$.

\ms

{\rm(ii)} The solution $\BPi(\cd)\in C([t,T];\dbR^{n\times 2n})$ to the ODE on $[t,T]$
\bel{16June29-Pi}\left\{\2n\ba{ll}
\ds\dot\BPi+\BPi A+\BA^\top\BPi+\BC^\top\BPi C+\BQ\BI_n
+\big(\BPi B+\BC^\top\BPi D+\BS^\top\BI_m\big)\Th=0,\\
\ns\ds\BPi(T)=\BG\BI_n,\ea\right.\ee
satisfies
\bel{16June29-21:40-1}\big[\BJ^\top(\BR\BI_m+\BD^\top\BPi D)\big]\Th
+\BJ^\top\big(\BB^\top\BPi+\BD^\top\BPi C+\BS\BI_n\big)=0,\ee
and the adapted solution $(\BBeta(\cd),\BBzeta(\cd))$ to the BSDE on $[t,T]$
\bel{16June29-BBeta}\left\{\2n\ba{ll}
\ds d\BBeta=-\big[\BA^\top\BBeta+\BC^\top\BBzeta+\big(\BPi B+\BC^\top\BPi D+\BS^\top\BI_m\big)v
+\BC^\top\BPi\si+\BPi b+\Bq\big]ds+\BBzeta dW,\\
\ns\ds\BBeta(T)=\Bg,\ea\right.\ee
satisfies
\bel{16June29-21:40-2}\big[\BJ^\top(\BR\BI_m+\BD^\top\BPi D)\big]v
+\BJ^\top\big(\BB^\top\BBeta+\BD^\top\BBzeta+\BD^\top\BPi\si+\Brho\big)=0.\ee
\et

\it Proof. \rm For any $x\in\dbR^n$, let $X(\cd)$, $\BPi(\cd)$, and $(\BBeta(\cd),\BBzeta(\cd))$
be the solutions to \rf{16June29-X}, \rf{16June29-Pi}, and \rf{16June29-BBeta}, respectively.
Let $u(\cd)$ be defined by \rf{16June29-17:40} and set
$$\BY=\BPi X+\BBeta,\qq \BZ=\BPi(C+D\Th)X+\BPi Dv+\BPi\si+\BBzeta.$$
Then $\BY(T)=\BG\BI_nX(T)+\Bg$, and
$$\ba{lll} 
\ds d\BY\4n&=&\4n\ds \dot\BPi Xds+\BPi dX+d\BBeta\\
\ns\4n&=&\4n\ds \big[\dot\BPi X+\BPi(A+B\Th)X+\BPi Bv+\BPi b-\BA^\top\BBeta-\BC^\top\BBzeta\\
\ns\4n&~&\4n\ds\qq~-\big(\BPi B+\BC^\top\BPi D+\BS^\top\BI_m\big)v-\BC^\top\BPi\si-\BPi b-\Bq\big]ds\\
\ns\4n&~&\4n\ds+\,\big[\BPi(C+D\Th)X+\BPi Dv+\BPi\si+\BBzeta\big]dW\\
\ns\4n&=&\4n\ds \big[-\big(\BA^\top\BPi+\BC^\top\BPi C+\BQ\BI_n+\BC^\top\BPi D\Th+\BS^\top\BI_m\Th\big)X
-\BA^\top\BBeta-\BC^\top\BBzeta\\
\ns\4n&~&\4n\ds\qq~-\big(\BC^\top\BPi D+\BS^\top\BI_m\big)v-\BC^\top\BPi\si-\Bq\big]ds+\BZ dW\\
\ns\4n&=&\4n\ds \big\{-\BA^\top(\BPi X+\BBeta)-\BQ\BI_n X-\BC^\top\big[\BPi(C+D\Th)X+\BPi Dv+\BPi\si+\BBzeta\big]\\
\ns\4n&~&\4n\ds\qq~-\BS^\top\BI_m(\Th X+v)-\Bq\big\}ds+\BZ dW\\
\ns\4n&=&\4n\ds \big(-\BA^\top\BY-\BQ\BI_n X-\BC^\top\BZ-\BS^\top\BI_m u-\Bq\big)ds+\BZ dW. \ea$$
This shows that $(X(\cd),\BY(\cd),\BZ(\cd),u(\cd))$ satisfies the FBSDE \rf{FBSDE**}.
According to Theorem \ref{Theorem 4.1}, the process $u(\cd)$ defined by \rf{16June29-17:40} is an
open-loop Nash equilibrium for $(t,x)$ if and only if (i) holds and
$$\ba{lll}
\ds0\4n&=&\4n\ds \BJ^\top\big(\BB^\top\BY+\BD^\top\BZ+\BS\BI_nX+\BR\BI_m u+\Brho\big)\\
\ns\4n&=&\4n\ds \BJ^\top\big\{\BB^\top(\BPi X+\BBeta)+\BD^\top[\BPi(C+D\Th)X+\BPi Dv+\BPi\si+\BBzeta]+\BS\BI_nX+\BR\BI_m(\Th X+v)+\Brho\big\}\\
\ns\4n&=&\4n\ds \BJ^\top\big[\BB^\top\BPi+\BD^\top\BPi C+\BS\BI_n+(\BR\BI_m+\BD^\top\BPi D)\Th\big]X\\
\ns\4n&~&\4n\ds\qq~ +\BJ^\top\big[\BB^\top\BBeta+\BD^\top\BBzeta+\BD^\top\BPi\si+\Brho+(\BR\BI_m+\BD^\top\BPi D)v\big].
\ea$$
Since the initial state $x$ is arbitrary and
$\BJ^\top[\BB^\top\BBeta+\BD^\top\BBzeta+\BD^\top\BPi\si+\Brho+(\BR\BI_m+\BD^\top\BPi D)v]$
is independent of $x$, the above leads to \rf{16June29-21:40-1} and \rf{16June29-21:40-2}.
\endpf

\ms

Let us write \rf{16June29-Pi}--\rf{16June29-21:40-2} componentwise as follows: For $i=1,2$,
\begin{eqnarray}
&\label{16June30-Pi}\left\{\2n\ba{ll}
\ds\dot\Pi_i+\Pi_iA+A^\top\Pi_i+C^\top\Pi_iC+Q^i+\big[\Pi_iB+C^\top\Pi_iD+(S^i)^\top\big]\Th=0,\\
\ns\ds\Pi_i(T)=G^i,\ea\right.&\\
\ns&\label{16June30-Pi-yueshu}\begin{pmatrix}R_1^1+D_1^\top\Pi_1D \\ R_2^2+D_2^\top\Pi_2D\end{pmatrix}\Th
+\begin{pmatrix}B_1^\top\Pi_1+D_1^\top\Pi_1C+S_1^1 \\B_2^\top\Pi_2+D_2^\top\Pi_2C+S_2^2\end{pmatrix}=0,&\\
\ns&\label{16June30-eta}\left\{\2n\ba{ll}
\ds d\eta_i=-\Big\{A^\top\eta_i+C^\top\zeta_i+\big[\Pi_iB+C^\top\Pi_iD+(S^i)^\top\big]v
+C^\top\Pi_i\si+\Pi_ib+q^i\Big\}ds+\zeta_i dW,\\
\ns\ds\eta_i(T)=g^i,\ea\right.&\\
\ns&\label{16June30-eta-yueshu}\begin{pmatrix}R_1^1+D_1^\top\Pi_1D \\ R_2^2+D_2^\top\Pi_2D\end{pmatrix}v
+\begin{pmatrix}B_1^\top\eta_1+D_1^\top\zeta_1+D_1^\top\Pi_1\si+\rho_1^1\\
                B_2^\top\eta_2+D_2^\top\zeta_2+D_2^\top\Pi_2\si+\rho_2^2\end{pmatrix}=0.&
\end{eqnarray}
Noting the relation \rf{16June30-Pi-yueshu}, one sees the equations for $\Pi_1(\cd)$ and $\Pi_2(\cd)$
are coupled and none of them is symmetric. Consequently, $\Pi_1(\cd)$ and $\Pi_2(\cd)$ are not symmetric
in general. Whereas the Riccati equations \rf{Riccati-i} for $P_i(\cd)$ $(i=1,2)$ are symmetric.
This is the main reason that the closed-loop representation of open-loop Nash equilibria is different
from the outcome of closed-loop Nash equilibria.

\section{Zero-Sum Cases}

In the previous section, we have seen that for Problem (SDG), the closed-loop representation of
open-loop Nash equilibria is different from the outcome of closed-loop Nash equilibria in general.
Now we would like to take a look at the situation for LQ stochastic two-person zero-sum differential games.
In this case, Nash equilibria are usually called saddle points. According to \rf{16Apr2-1+2=0}, we have
\bel{16July1-16:00}\ba{ll}
\ds G^1=-G^2\equiv G, \q g^1=-g^2\equiv g, \q Q^1(\cd)=-Q^2(\cd)\equiv Q(\cd), \q q^1(\cd)=-q^2(\cd)\equiv q(\cd),\\
\ns\ds\begin{pmatrix}R^1_{11}(\cd)&R^1_{12}(\cd) \\ R^1_{21}(\cd)&R^1_{22}(\cd)\end{pmatrix}
\equiv-\begin{pmatrix}R^2_{11}(\cd)&R^2_{12}(\cd) \\ R^2_{21}(\cd)&R^2_{22}(\cd)\end{pmatrix}
\equiv\begin{pmatrix}R_{11}(\cd)&R_{12}(\cd) \\ R_{21}(\cd)&R_{22}(\cd)\end{pmatrix}
\equiv\begin{pmatrix}R_1(\cd) \\ R_2(\cd)\end{pmatrix}\equiv R(\cd),\\
\ns\ds\begin{pmatrix}S^1_1(\cd) \\ S^1_2(\cd)\end{pmatrix}=-\begin{pmatrix}S^2_1(\cd) \\ S^2_2(\cd)\end{pmatrix}
\equiv\begin{pmatrix}S_1(\cd) \\ S_2(\cd)\end{pmatrix}\equiv S(\cd),
\q\begin{pmatrix}\rho^1_1(\cd) \\ \rho^1_2(\cd)\end{pmatrix}
=-\begin{pmatrix}\rho^2_1(\cd) \\ \rho^2_2(\cd)\end{pmatrix}
\equiv\begin{pmatrix}\rho_1(\cd) \\ \rho_2(\cd)\end{pmatrix}\equiv\rho(\cd),
\ea\ee
and
$$\ba{ll}
\ds J^1(t,x;u_1(\cd),u_2(\cd))=-J^2(t,x;u_1(\cd),u_2(\cd))\\
\ns\ds=\dbE\Big\{\lan GX(T),X(T)\ran+2\lan g,X(T)\ran\\
\ns\ds\qq~+\int_t^T\[\lan
{\sc\begin{pmatrix}\sc Q(s)   &\1n\sc S_1(s)^\top &\1n\sc S_2(s)^\top\\
                   \sc S_1(s) &\1n\sc R_{11}(s)   &\1n\sc R_{12}(s)  \\
                   \sc S_2(s) &\1n\sc R_{21}(s)   &\1n\sc R_{22}(s)  \end{pmatrix}}
{\sc\begin{pmatrix}\sc X(s) \\ \sc u_1(s) \\ \sc u_2(s)\end{pmatrix}},
{\sc\begin{pmatrix}\sc X(s) \\ \sc u_1(s) \\ \sc u_2(s)\end{pmatrix}}\ran
+2\lan
{\sc\begin{pmatrix}\sc q(s) \\ \sc\rho_1(s) \\ \sc\rho_2(s)\end{pmatrix},
    \begin{pmatrix}\sc X(s) \\ \sc u_1(s)   \\ \sc u_2(s)  \end{pmatrix}}\ran\]ds\Big\}\\
\ns\ds\equiv J(t,x;u_1(\cd),u_2(\cd)). \ea$$
Let $(\Th(\cd),v(\cd))\in\cQ[t,T]\times\cU[t,T]$ and assume the open-loop saddle points of
Problem {\rm(SDG)} on $[t,T]$ admit the closed-loop representation \rf{16June29-17:40}.
The equations \rf{16June30-Pi} $(i=1,2)$ for $\Pi_1(\cd)$ and $\Pi_2(\cd)$ now become
$$\left\{\2n\ba{ll}
\ds\dot\Pi_1+\Pi_1 A+A^\top\Pi_1+C^\top\Pi_1 C+Q+\big(\Pi_1 B+C^\top\Pi_1 D+S^\top\big)\Th=0,\\
\ns\ds\Pi_1(T)=G,\ea\right.$$
and
$$\left\{\2n\ba{ll}
\ds\dot\Pi_2+\Pi_2 A+A^\top\Pi_2+C^\top\Pi_2 C-Q+\big(\Pi_2 B+C^\top\Pi_2 D-S^\top\big)\Th=0,\\
\ns\ds\Pi_2(T)=-G,\ea\right.$$
respectively. Obviously, both $\Pi_1(\cd)$ and $-\Pi_2(\cd)$ satisfy
\bel{16June30-17:00}\left\{\2n\ba{ll}
\ds\dot\Pi+\Pi A+A^\top\Pi+C^\top\Pi C+Q+\big(\Pi B+C^\top\Pi D+S^\top\big)\Th=0,\\
\ns\ds\Pi(T)=G.\ea\right.\ee
Thus, $\Pi_1(\cd)=-\Pi_2(\cd)\equiv\Pi(\cd)$, and \rf{16June30-Pi-yueshu} becomes
$$\begin{pmatrix}R_1+D_1^\top\Pi D \\ -R_2-D_2^\top\Pi D\end{pmatrix}\Th
+\begin{pmatrix}B_1^\top\Pi+D_1^\top\Pi C+S_1 \\-B_2^\top\Pi-D_2^\top\Pi C-S_2\end{pmatrix}=0,$$
or equivalently,
$$(R+D^\top\Pi D)\Th+B^\top\Pi+D^\top\Pi C+S=0.$$
This is also equivalent to
\bel{16June30-16:43}\left\{\2n\ba{ll}
\ds \sR(B^\top\Pi+D^\top\Pi C+S)\subseteq\sR(R+D^\top\Pi D),\qq\ae~s\in[t,T],\\
\ns\ds (R+D^\top\Pi D)^\dag(B^\top\Pi+D^\top\Pi C+S)\in L^2(t,T;\dbR^{m\times n}),
\ea\right.\ee
and
\bel{16June30-16:40}\Th=-(R+D^\top\Pi D)^\dag(B^\top\Pi+D^\top\Pi C+S)
+[I-(R+D^\top\Pi D)^\dag(R+D^\top\Pi D)]\th,\ee
for some $\th(\cd)\in L^2(t,T;\dbR^{m\times n})$. Upon substitution of \rf{16June30-16:40}
into \rf{16June30-17:00}, the latter becomes
\bel{16June30-20:00}\left\{\2n\ba{ll}
\ds\dot\Pi+\Pi A+A^\top\Pi+C^\top\Pi C+Q\\
\ns\ds\q-\,(\Pi B+C^\top\Pi D+S^\top)(R+D^\top\Pi D)^\dag(B^\top\Pi+D^\top\Pi C+S)=0,\q~ s\in[t,T],\\
\ns\ds\Pi(T)=G,\ea\right.\ee
with constraints \rf{16June30-16:43}. Note that equation \rf{16June30-20:00} is symmetric. Likewise,
we have $(\eta_1(\cd),\z_1(\cd))=-(\eta_2(\cd),\z_2(\cd))\equiv(\eta_\Pi(\cd),\z_\Pi(\cd))$ satisfying
\bel{16June30-20:14}\left\{\2n\ba{ll}
\ds d\eta_\Pi=-\Big\{\big[A^\top\1n-(\Pi B+C^\top\Pi D+S^\top)(R+D^\top\Pi D)^\dag B^\top\big]\eta_\Pi\\
\ns\ds\qq\qq~+\big[C^\top\1n-(\Pi B+C^\top\Pi D+S^\top)(R+D^\top\Pi D)^\dag D^\top\big]\z_\Pi\\
\ns\ds\qq\qq~+\big[C^\top\1n-(\Pi B+C^\top\Pi D+S^\top)(R+D^\top\Pi D)^\dag D^\top\big]\Pi\si\\
\ns\ds\qq\qq~-(\Pi B+C^\top\Pi D+S^\top)(R+D^\top\Pi D)^\dag\rho+\Pi b+q\Big\}ds+\z_\Pi dW,\q~s\in[t,T],\\
\ns\ds\eta_\Pi(T)=g,\ea\right.\ee
with constraints
\bel{16June30-20:16}\left\{\2n\ba{ll}
\ds B^\top\eta_\Pi+D^\top\z_\Pi+D^\top\Pi\si+\rho\in\sR(R+D^\top\Pi D),\qq\ae~s\in[t,T],~\as\\
\ns\ds (R+D^\top\Pi D)^\dag(B^\top\eta_\Pi+D^\top\z_\Pi+D^\top\Pi\si+\rho)\in L_\dbF^2(t,T;\dbR^m),
\ea\right.\ee
and in this case,
$$v=-(R+D^\top\Pi D)^\dag(B^\top\eta_\Pi+D^\top\z_\Pi+D^\top\Pi\si+\rho)
+\big[I-(R+D^\top\Pi D)^\dag(R+D^\top\Pi D)\big]\n,$$
for some $\n(\cd)\in L_\dbF^2(t,T;\dbR^m)$.
To summarize, we have the following result for LQ stochastic two-person zero-sum differential games.

\bt{bt-16June30-18:00}\sl Let {\rm(G1)--(G2)} and \rf{16July1-16:00} hold.
Then the open-loop saddle points of Problem {\rm(SDG)} on $[t,T]$ admit a closed-loop representation
if and only if the following hold:

\ms

{\rm(i)} The following convexity-concavity condition holds: For $i=1,2$,
\bel{16July1-17:28}\ba{ll}
\ds(-1)^{i-1}\dbE\Big\{\int_t^T\[\lan Q(s)X_i(s),X_i(s)\ran+2\lan S_i(s)X_i(s),u_i(s)\ran
+\lan R_{ii}(s)u_i(s),u_i(s)\ran\]ds\\
\ns\ds\qq\qq\qq\qq+\,\lan GX_i(T),X_i(T)\ran\Big\}\ges0,\qq\forall\,u_i(\cd)\in\cU_i[t,T],\ea\ee
where $X_i(\cd)$ is the solution to FSDE \rf{homogeneous}.

\ms

{\rm(ii)} The Riccati equation \rf{16June30-20:00} admits a solution $\Pi(\cd)\in C([t,T];\dbS^n)$ such that
\rf{16June30-16:43} holds, and the adapted solution of \rf{16June30-20:14} satisfies \rf{16June30-20:16}.

\ms

In the above case, all the closed-loop representations of open-loop saddle points are given by
$$\ba{lll}
\ds u\4n&=&\4n\ds \Big\{\1n-(R+D^\top\Pi D)^\dag(B^\top\Pi+D^\top\Pi C+S)
+\big[I-(R+D^\top\Pi D)^\dag(R+D^\top\Pi D)\big]\th\Big\}X\\
\ns\4n&~&\4n\ds-\,(R+D^\top\Pi D)^\dag(B^\top\eta_\Pi+D^\top\z_\Pi+D^\top\Pi\si+\rho)
+\big[I-(R+D^\top\Pi D)^\dag(R+D^\top\Pi D)\big]\n,
\ea$$
where $\th(\cd)\in L^2(t,T;\dbR^{m\times n})$ and $\n(\cd)\in L_\dbF^2(t,T;\dbR^m)$.
\et

\it Proof. \rm The result can be proved by combining Theorem \ref{bt-16June29-17:30} and the previous argument.
We leave the details to the interested reader.
\endpf

\ms

Now let us recall from \cite{Sun-Yong 2014} the characterization of closed-loop saddle points
of LQ stochastic two-person zero-sum differential games.

\bt{bt-16July1-15:00} \sl Let {\rm(G1)--(G2)} and \rf{16July1-16:00} hold. Then Problem {\rm(SDG)}
admits a closed-loop saddle point on $[t,T]$ if and only if the following hold:

\ms

{\rm(i)} The Riccati equation
\bel{16July1-Ric}\left\{\2n\ba{ll}
\ds \dot P+PA+A^\top P+C^\top PC+Q\\
\ns\ds\q-\,(PB+C^\top PD+S^\top)(R+D^\top PD)^\dag(B^\top P+D^\top PC+S)=0,\q~s\in[t,T],\\
\ns\ds P(T)=G,\ea\right.\ee
admits a solution $P(\cd)\in C([t,T];\dbS^n)$ such that the following hold:
\begin{eqnarray}
&\label{16July1-yushu*}\left\{\2n\ba{ll}
\ds \sR(B^\top P+D^\top PC+S)\subseteq\sR(R+D^\top PD),\qq\ae~s\in[t,T],\\
\ns\ds (R+D^\top PD)^\dag(B^\top P+D^\top PC+S)\in L^2(t,T;\dbR^{m\times n}),
\ea\right.&\\
\ns&\label{16July1-yushu**}R_{11}+D_1^\top PD_1\ges0,\qq R_{22}+D_2^\top PD_2\les0,\qq\ae~s\in[t,T].&
\end{eqnarray}

{\rm(ii)} The adapted solution $(\eta_P(\cd),\z_P(\cd))$ of the BSDE on $[t,T]$
\bel{16July2-eta_P}\left\{\2n\ba{ll}
\ds d\eta_P=-\Big\{\big[A^\top\1n-(PB+C^\top PD+S^\top)(R+D^\top PD)^\dag B^\top\big]\eta_P\\
\ns\ds\qq\qq~+\big[C^\top\1n-(PB+C^\top PD+S^\top)(R+D^\top PD)^\dag D^\top\big]\z_P\\
\ns\ds\qq\qq~+\big[C^\top\1n-(PB+C^\top PD+S^\top)(R+D^\top PD)^\dag D^\top\big]P\si\\
\ns\ds\qq\qq~-(PB+C^\top PD+S^\top)(R+D^\top PD)^\dag\rho+Pb+q\Big\}ds+\z_P dW,\\
\ns\ds\eta_P(T)=g,\ea\right.\ee
satisfies
\bel{}\left\{\2n\ba{ll}
\ds B^\top\eta_P+D^\top\z_P+D^\top P\si+\rho\in\sR(R+D^\top PD),\qq\ae~s\in[t,T],~\as\\
\ns\ds (R+D^\top PD)^\dag(B^\top\eta_P+D^\top\z_P+D^\top P\si+\rho)\in L_\dbF^2(t,T;\dbR^m).
\ea\right.\ee

In this case, the closed-loop saddle point $(\Th^*(\cd),v^*(\cd))$ admits the following representation:
\bel{16July2-15:30}\left\{\2n\ba{cll}
\ds\Th^*\4n&=&\4n\ds -(R+D^\top PD)^\dag(B^\top P+D^\top PC+S)
+\big[I-(R+D^\top PD)^\dag(R+D^\top PD)\big]\th,\\
\ns v^*\4n&=&\4n\ds -(R+D^\top PD)^\dag(B^\top\eta_P+D^\top\z_P+D^\top P\si+\rho)
+\big[I-(R+D^\top PD)^\dag(R+D^\top PD)\big]\n,
\ea\right.\ee
where $\th(\cd)\in L^2(t,T;\dbR^{m\times n})$ and $\n(\cd)\in L_\dbF^2(t,T;\dbR^m)$.
\et

Comparing Theorems \ref{bt-16June30-18:00} and \ref{bt-16July1-15:00}, one may ask:
For LQ stochastic two-person zero-sum differential games, when both the closed-loop representation of
open-loop saddle points and the closed-loop saddle point exist, does the closed-loop representation
coincide with the outcome of the closed-loop saddle point?
The answer to this question is affirmative, as shown by the following result.

\bt{bt-16July1-17:00}\sl Let {\rm(G1)--(G2)} and \rf{16July1-16:00} hold. If both the closed-loop
representation of open-loop saddle points and the closed-loop saddle point exist on $[t,T]$,
then the closed-loop representation coincides with the outcome of the closed-loop saddle point.
\et

\it Proof. \rm The proof is immediate from Theorems \ref{bt-16June30-18:00} and \ref{bt-16July1-15:00},
once we show that the solution $\Pi(\cd)$ to the Riccati equation \rf{16June30-20:00} with constraints
\rf{16June30-16:43} coincides with the solution $P(\cd)$ to \rf{16July1-Ric} with constraints
\rf{16July1-yushu*}--\rf{16July1-yushu**}.

\ms

First, we note that if the convexity-concavity condition \rf{16July1-17:28} holds for initial time $t$,
it also holds for any $t^\prime\in[t,T]$. Indeed, for any $t^\prime\in[t,T]$,
and any $u_1(\cd)\in\cU_1[t^\prime,T]$, let $X_1(\cd)$ be the solution to
$$\left\{\2n\ba{ll}
\ds dX_1(s)=\big[A(s)X_1(s)+B_1(s)u_1(s)\big]ds+\big[C(s)X_1(s)+D_1(s)u_1(s)\big]dW(s),\q~s\in[t^\prime,T], \\
\ns\ds X_1(t^\prime)=0,\ea\right.$$
and define the {\it zero-extension}
of $u_1(\cd)$ as follows:
$$[\,0I_{[t,t^\prime)}\oplus u_1](s)=\left\{\2n\ba{ll}0,& s\in[t,t^\prime),\\
\ns\ds u_1(s),& s\in[t^\prime,T].\ea\right.$$
Then $\tilde{u}_1(\cd)\equiv[\,0I_{[t,t^\prime)}\oplus u_1](\cd)\in\cU_1[t,T]$, and due to the initial state being 0,
the solution $\wt X_1(s)$ of
$$\left\{\2n\ba{ll}
\ds d\wt X_1(s)=\big[A(s)\wt X_1(s)+B_1(s)\tilde{u}_1(s)\big]ds+\big[C(s)\wt X_1(s)+D_1(s)\tilde{u}_1(s)\big]dW(s),\q s\in[t,T],\\
\ns\ds \wt X_1(t)=0,\ea\right.$$
satisfies
$$\wt X_1(s)=\left\{\2n\ba{ll}0,& s\in[t,t^\prime),\\
\ns\ds X_1(s),& s\in[t^\prime,T].\ea\right.$$
Hence,
$$\ba{ll}
\ds\dbE\Big\{\int_{t^\prime}^T\[\lan QX_1,X_1\ran+2\lan S_1X_1,u_1\ran+\lan R_{11}u_1,u_1\ran\]ds
+\lan GX_1(T),X_1(T)\ran\Big\}\\
\ns\ds=\dbE\Big\{\int_t^T\[\blan Q\wt X_1,\wt X_1\bran+2\blan S_1\wt X_1,\tilde{u}_1\bran
+\blan R_{11}\tilde{u}_1,\tilde{u}_1\bran\]ds+\blan G\wt X_1(T),\wt X_1(T)\bran\Big\}\ges0.
\ea$$
This proves the case $i=1$. The case $i=2$ can be treated similarly.

\ss

Now let $(\Th^*(\cd),v^*(\cd))$ be a closed-loop saddle point of Problem (SDG) on $[t,T]$.
Under the assumption of the theorem, it is clear from Theorem \ref{bt-16June30-18:00} that
for any initial pair $(t^\prime,x)$ with $t^\prime\in[t,T]$, the outcome
$$u^*(s)=\Th^*(s)X^*(s)+v^*(s),\qq s\in[t^\prime,T]$$
of $(\Th^*(\cd),v^*(\cd))$ is an open-loop saddle point for $(t^\prime,x)$, where $X^*(\cd)$
is the solution to
$$\left\{\2n\ba{ll}
\ds dX^*(s)=\big\{[A(s)+B(s)\Th^*(s)]X^*(s)+B(s)v^*(s)+b(s)\big\}ds\\
\ns\ds\qq\qq~+\big\{[C(s)+D(s)\Th^*(s)]X^*(s)+D(s)v^*(s)+\si(s)\big\}dW(s),\qq s\in[t^\prime,T],\\
\ns\ds X^*(t^\prime)= x. \ea\right.$$
By Theorem \ref{bt-16July1-15:00}, $(\Th^*(\cd),v^*(\cd))$ admits the representation \rf{16July2-15:30},
and a straightforward calculation shows that
$$\ba{ll}
\ds\dot P+P(A+B\Th^*)+(A+B\Th^*)^\top\1n P+(C+D\Th^*)^\top\1n P(C+D\Th^*)\\
\ns\ds\q+\,(\Th^*)^\top R\Th^*+S^\top\Th^*+(\Th^*)^\top S+Q=0,\ea$$
and that the adapted solution $(\eta_P(\cd),\z_P(\cd))$ of \rf{16July2-eta_P} satisfies
$$ d\eta_P=-\big[(A+B\Th^*)^\top\eta_P+(C+D\Th^*)^\top\z_P+(C+D\Th^*)^\top P\si
+(\Th^*)^\top\rho+Pb+q\big]ds+\z_P dW.$$
Then applying It\^o's formula to $s\mapsto\lan P(s)X^*(s),X^*(s)\ran+2\lan \eta_P(s),X^*(s)\ran$
and noting that
$$(R+D^\top PD)\Th^*+B^\top P+D^\top PC+S=0,$$
we have
\bel{16July2-17:44}\ba{ll}
\ds J(t^\prime,x;u^*(\cd))=J(t^\prime,x;\Th^*(\cd)X^*(\cd)+v^*(\cd))\\
\ns\ds=\dbE\Big\{\lan GX^*(T),X^*(T)\ran+2\lan g,X^*(T)\ran+\int_{t^\prime}^T\[\lan QX^*,X^*\ran+2\lan SX^*,\Th^*X^*+v^*\ran\\
\ns\ds\qq~+\lan R(\Th^*X^*+v^*),\Th^*X^*+v^*\ran+2\lan q,X^*\ran+2\lan\rho,\Th^*X^*+v^*\ran\]ds\Big\}\\
\ns\ds=\dbE\Big\{\lan P(t^\prime)x,x\ran+2\lan\eta_P(t^\prime),x\ran+\int_{t^\prime}^T\[\lan\dot PX^*,X^*\ran+2\lan PX^*,(A+B\Th^*)X^*+Bv^*+b\ran\\
\ns\ds\qq~+\lan P[(C+D\Th^*)X^*+Dv^*+\si],(C+D\Th^*)X^*+Dv^*+\si\ran\\
\ns\ds\qq~-2\blan(A+B\Th^*)^\top\eta_P+(C+D\Th^*)^\top\z_P+(C+D\Th^*)^\top P\si+(\Th^*)^\top\rho+Pb+q,X^*\bran\\
\ns\ds\qq~+2\lan\eta_P,(A+B\Th^*)X^*+Bv^*+b\ran+2\lan\z_P,(C+D\Th^*)X^*+Dv^*+\si\ran\\
\ns\ds\qq~+\blan\big[Q+S^\top\Th^*+(\Th^*)^\top S+(\Th^*)^\top R\Th^*\big]X^*,X^*\bran+2\lan (R\Th^*+S)X^*,v^*\ran\\
\ns\ds\qq~+2\blan q+(\Th^*)^\top\rho,X^*\bran+\lan Rv^*,v^*\ran+2\lan\rho,v^*\ran\]ds\Big\}\\
\ns\ds=\dbE\Big\{\lan P(t^\prime)x,x\ran+2\lan\eta_P(t^\prime),x\ran+\int_{t^\prime}^T\[\lan P\si,\si\ran+2\lan\eta_P,b\ran+2\lan\z_P,\si\ran\\
\ns\ds\qq~+\blan(R+D^\top PD)v^*,v^*\bran+2\blan B^\top\eta_P+D^\top\z_P+D^\top P\si+\rho,v^*\bran\]ds\Big\}.
\ea\ee

\ss

Next, let $\th(\cd)\in L^2(t,T;\dbR^{m\times n}), \n(\cd)\in L_\dbF^2(t,T;\dbR^m)$ and denote
$$\left\{\2n\ba{cll}
\ds \Th\4n&=&\4n\ds -(R+D^\top\Pi D)^\dag(B^\top\Pi+D^\top\Pi C+S)
+[I-(R+D^\top\Pi D)^\dag(R+D^\top\Pi D)]\th,\\
\ns v\4n&=&\4n\ds -(R+D^\top\Pi D)^\dag(B^\top\eta_\Pi+D^\top\z_\Pi+D^\top\Pi\si+\rho)
+[I-(R+D^\top\Pi D)^\dag(R+D^\top\Pi D)]\n.
\ea\right.$$
For any initial pair $(t^\prime,x)$ with $t^\prime\in[t,T]$, define $u(\cd)\in\cU[t^\prime,T]$ by
$$u(s)=\Th(s)X(s)+v(s),\qq s\in[t^\prime,T],$$
with $X(\cd)$ being the solution to
$$\left\{\2n\ba{ll}
\ds dX(s)=\big\{[A(s)+B(s)\Th(s)]X(s)+B(s)v(s)+b(s)\big\}ds\\
\ns\ds\qq\qq~+\big\{[C(s)+D(s)\Th(s)]X(s)+D(s)v(s)+\si(s)\big\}dW(s),\qq s\in[t^\prime,T],\\
\ns\ds X(t^\prime)= x. \ea\right.$$
By Theorem \ref{bt-16June30-18:00}, $u(\cd)$ is an open-loop saddle point for $(t^\prime,x)$,
and by a computation similar to \rf{16July2-17:44}, we obtain
\bel{16July2-18:15}\ba{ll}
\ds J(t^\prime,x;u(\cd))=\dbE\Big\{\lan\Pi(t^\prime)x,x\ran+2\lan\eta_\Pi(t^\prime),x\ran
+\int_{t^\prime}^T\[\lan\Pi\si,\si\ran+2\lan\eta_\Pi,b\ran+2\lan\z_\Pi,\si\ran\\
\ns\ds\qq\qq\qq\q~~+\blan(R+D^\top\Pi D)v,v\bran+2\blan B^\top\eta_\Pi+D^\top\z_\Pi+D^\top\Pi\si+\rho,v\bran\]ds\Big\}.
\ea\ee
Since both $u^*(\cd)\equiv(u_1^*(\cd)^\top,u_2^*(\cd)^\top)^\top$ and $u(\cd)\equiv(u_1(\cd)^\top,u_2(\cd)^\top)^\top$
are open-loop saddle points for $(t^\prime,x)$, we have
$$\ba{lll}
\ds J(t^\prime,x;u_1^*(\cd),u_2^*(\cd))\4n&\les&\4n\ds J(t^\prime,x;u_1(\cd),u^*_2(\cd))\les J(t^\prime,x;u_1(\cd),u_2(\cd))\\
\ns\4n&\les&\4n\ds J(t^\prime,x;u^*_1(\cd),u_2(\cd))\les J(t^\prime,x;u_1^*(\cd),u_2^*(\cd)).\ea$$
Therefore, $J(t^\prime,x;u^*(\cd))=J(t^\prime,x;u(\cd))$ for all $(t^\prime,x)$ with $t^\prime\in[t,T]$,
which, together with \rf{16July2-17:44} and \rf{16July2-18:15}, yields $\Pi(\cd)=P(\cd)$.
\endpf

\ms

Finally, we have the following corollary for Problem (SLQ), which should be but has not been stated in \cite{Sun-Li-Yong 2016}.

\bc{bc-16July11-15:30}\sl For Problem {\rm(SLQ)}, if the open-loop optimal controls admit
a closed-loop representation, then every open-loop optimal control must be an outcome of
a closed-loop optimal strategy.
\ec

\end{document}